\documentclass[preprint,12pt]{elsarticle}

\usepackage[T1]{fontenc}
\usepackage[latin1,utf8]{inputenc}
\usepackage{caption}
\usepackage{graphicx,url}
\usepackage{amsmath,amssymb,amsopn}
\usepackage{tabularx}
\usepackage{color}
\usepackage{calc}
\usepackage{longtable}
\usepackage{multirow}
\usepackage{mathrsfs}
\usepackage{array}
\usepackage{hyperref}
\usepackage{tikz}
\usepackage{booktabs}
\usepackage{pgfplotstable}
\usepackage{pgfplots}
\pgfplotsset{compat=1.18}
\usepackage[outline]{contour}
\usepackage[colorinlistoftodos,prependcaption]{todonotes}
\usepackage{calrsfs}
\usepackage{graphicx}
\usepgfplotslibrary{groupplots}
\usetikzlibrary{matrix}
\usetikzlibrary{automata, external, shapes.geometric, fit}
\usetikzlibrary{arrows.meta}
\usepackage{listings}
\usepackage{xcolor}

\usepackage[font=normalsize,labelfont=bf]{caption}

\lstdefinestyle{mystyle}{
  basicstyle=\ttfamily\small,
  backgroundcolor=\color{gray!10}, 
  frame=single,              
  breaklines=true,                  
  columns=fullflexible,
  keepspaces=true,  
}
\def\dv{\vec{d}}

\def\pv{\vec{p}}

\def\uv{\vec{u}}

\def\xv{\vec{x}}
\def\yv{\vec{y}}

\def\thetav{\vec{\theta}}
\def\xiv{\vec{\xi}}

\def\etav{\vec{\eta}}

\def\Hv{\vec{H}}

\def\Rv{\vec{R}}

\def\Hm{\mat{H}}
\def\Im{\mat{I}}

\def\R{\mathbb{R}}			
\def\N {\mathbb{N}}			
\def\E {\mathbb{E}}			

\def\D{\pazocal{D}}			
\def\L{\pazocal{L}}			%
\def\D{\pazocal{D}} 			
\def\F{\pazocal{F}}
\def\A{\pazocal{A}}
\def\B{\pazocal{B}}
\def\Y{\pazocal{Y}}
\def\V{\pazocal{V}}
\def\G{\pazocal{G}}
\def\K{\pazocal{K}}

\DeclareMathAlphabet{\pazocal}{OMS}{zplm}{m}{n}

\renewcommand{\vec}[1]{\boldsymbol{#1}}
\newcommand{\mat}[1]{\boldsymbol{\mathrm{#1}}}

\definecolor{myblack}{RGB}{53, 53, 53}
\definecolor{myblue}{RGB}{40, 75, 99}
\definecolor{myred}{RGB}{192, 50, 33}
\definecolor{myyellow}{RGB}{255, 166, 48}
\definecolor{mywhite}{RGB}{240, 237, 238}
\definecolor{mygreen}{RGB}{0, 102, 0}

\definecolor{green1}{RGB}{9, 82, 86}
\definecolor{green2}{RGB}{8, 127, 140}
\definecolor{green3}{RGB}{6, 167, 125}
\definecolor{green4}{RGB}{79, 109, 122}
\definecolor{green5}{RGB}{192, 214, 223}
\definecolor{violet}{RGB}{26,69,131}

\definecolor{checkgreen}{rgb}{0,0.6,0}
\definecolor{phase1}{rgb}{0.008,0.655,1.000}
\definecolor{phase2}{rgb}{0.016,0.75,0.700}
\definecolor{phase3}{rgb}{0.929,0.35,0.700}
\definecolor{icsyellow}{cmyk}{0.00,0.11,0.53,0.00}

\usepackage{color, colortbl}
\definecolor{Gray}{gray}{0.9}

\definecolor{checkgreen}{rgb}{0,0.6,0}
\definecolor{phase1}{rgb}{0.008,0.655,1.000}
\definecolor{phase2}{rgb}{0.016,0.75,0.700}
\definecolor{phase3}{rgb}{0.929,0.35,0.700}
\definecolor{icsyellow}{cmyk}{0.00,0.11,0.53,0.00}
\definecolor{green1}{RGB}{9, 82, 86}
\definecolor{green2}{RGB}{8, 127, 140}
\definecolor{green3}{RGB}{6, 167, 125}
\definecolor{green4}{RGB}{79, 109, 122}
\definecolor{green5}{RGB}{192, 214, 223}
\definecolor{violet}{RGB}{26,69,131}
\definecolor{mypurple}{RGB}{150, 0, 180}
\definecolor{olive}{HTML}{808000}

\usepackage[noabbrev]{cleveref}
\usepackage{lineno}

\usepackage{lipsum}
\usepackage{amsfonts}
\usepackage{graphicx}
\usepackage{epstopdf}
\usepackage{algorithm}
\usepackage{algorithmic}
\usepackage{enumitem}

\usepackage{amsopn}

\usepackage{amsthm}
\newtheorem{theorem}{Theorem}
\theoremstyle{remark}
\newtheorem{remark}[theorem]{Remark}

\usetikzlibrary{intersections}
\usepgfplotslibrary{fillbetween}
\usepackage{multirow}
\usepackage{stmaryrd}
\usepackage{xfrac}
\usepackage{caption}
\usepackage{subcaption}
\usepackage{lmodern}
\usetikzlibrary{decorations.text}
\usetikzlibrary{patterns}
\usetikzlibrary{patterns.meta}
\usepackage[normalem]{ulem}

\definecolor{blackmy}{RGB}{38, 70, 83}
\definecolor{greenmy}{RGB}{42, 167, 143}
\definecolor{yellowmy}{RGB}{233, 196, 106}
\definecolor{redmy}{RGB}{253,127,111}

\definecolor{bluemy}{RGB}{40, 75, 99}
\definecolor{myred}{RGB}{192, 50, 33}
\definecolor{brownmy}{RGB}{244, 162, 97}

\definecolor{blue1}{RGB}{1, 58, 99}
\definecolor{blue2}{RGB}{70, 143, 175}
\definecolor{blue3}{RGB}{137, 194, 217}

\definecolor{red1}{RGB}{176, 71, 89}
\definecolor{red3}{RGB}{249, 155, 125}
\definecolor{red2}{RGB}{231, 97, 97}

\definecolor{green1}{RGB}{112, 151, 117}
\definecolor{green2}{RGB}{143, 185, 150}
\definecolor{green3}{RGB}{161, 204, 165}

\DeclareMathOperator*{\argmin}{arg\,min}

\newcommand{\algorithmiccommentMine}[1]{\bgroup\hfill$\triangleright$~{#1}\egroup}

\newcommand\newtext[1]{#1}

\newcommand\oldtext[1]{}
\newcommand\cancel[1]{}

\journal{Computer Methods in Applied Mechanics and Engineering}

\begin{document} 

\begin{frontmatter}



	\title{Two-level overlapping Additive Schwarz preconditioner for training Scientific Machine learning Applications}


	\author[brown]{Youngkyu Lee} 
	\author[brown,INP,ANITI]{Alena Kopani\v{c}\'akov\'a} 
	\author[brown]{George Em Karniadakis~\fnref{2}} 

	\affiliation[brown]{organization={Division of Applied Mathematics, Brown University},
		addressline={170 Hope St},
		city={Providence},
		postcode={02906},
		state={RI},
		country={USA}}
	\affiliation[INP]{organization={Institut national polytechnique de Toulouse/Enseeiht - Institut de Recherche en Informatique de Toulouse},
		addressline={2 Rue Charles Camichel},
		city={Toulouse},
		postcode={31000},
		country={France}}
	\affiliation[ANITI]{organization={Artificial and Natural Intelligence Toulouse Institute (ANITI)},
		addressline={3 Rue Tarfaya},
		city={Toulouse},
		postcode={314000},
		country={France}}        
	\fntext[2]{Corresponding author: george\_karniadakis@brown.edu (George Em Karniadakis).}

	\begin{abstract}
		We introduce a novel two-level overlapping additive Schwarz preconditioner for accelerating the training of scientific machine learning applications.
		The design of the proposed preconditioner is motivated by the nonlinear two-level overlapping additive Schwarz preconditioner.
		The neural network parameters are decomposed into groups (subdomains) with overlapping regions.
		In addition, the network's feed-forward structure is indirectly imposed through a novel subdomain-wise synchronization strategy and a coarse-level training step.
		Through a series of numerical experiments, which consider physics-informed neural networks and operator learning approaches, we demonstrate that the proposed two-level preconditioner significantly speeds up the convergence of the standard (LBFGS) optimizer while also yielding more accurate machine learning models.
		Moreover, the devised preconditioner is designed to take advantage of model-parallel computations, which can further reduce the training time.
	\end{abstract}




\end{frontmatter}



\section{Introduction}
Deep neural networks (DNNs) are universal approximators, capable of approximating any continuous function defined on a compact domain to arbitrary accuracy~\cite{barron1993universal,chen1990constructive,siegel2020approximation}.
Therefore, they have been recently widely used to obtain ansatz spaces for the solutions of partial differential equations (PDEs) for a wide range of scientific and engineering problems; see~\cite{karniadakis2021physics} for an overview.
The popularity of these DNN approaches can be largely attributed to their simplicity, mesh-free nature, and ability to incorporate data~\cite{raissi2019physics, sirignano2018dgm}, enabling them to effectively tackle complex, high-dimensional, forward and inverse problems in single and multi-query scenarios~\cite{mishra2022estimates}.

Two prominent approaches for solving parametric PDEs using DNNs appear in the literature.
The first approach approximates a solution of a given PDE for a particular choice of parameters using a DNN~\cite{raissi2019physics,sirignano2018dgm}.
Popular representatives of this approach are physics informed neural networks (PINNs)~\cite{raissi2019physics}, which train the DNNs by minimizing the mean square error that incorporates the PDE residual, boundary/initial conditions, and/or observed data.
The second approach considers learning the differential operator~\cite{lu2021learning,wang2021learning} from the parametrized source terms, boundary/initial conditions, or material properties.
The operator learning DNNs, such as DeepONet (DON)~\cite{lu2021learning} or FNO~\cite{li2021fourier}, are typically trained by minimizing the mean square error between the network's output and the solution generated using a high-fidelity method.

Both, physics informed and operator learning, approaches are commonly trained using gradient-based optimizers, e.g., Adam~\cite{kingma2014adam} or LBFGS~\cite{liu1989limited}.
However, the computational cost of these optimizers can be exorbitant due to the non-convexity and ill-conditioning of the underlying optimization problem.
To address this problem, there has been research on improving these optimizers~\cite{Liu2020On,loshchilov2018decoupled,Al-Baali03092014}.
This is especially the case when a highly accurate solution is required.
In order to reduce the training time, the parallelization of the DNN training process has gained significant attention in the literature; see, for example~\cite{ben2019demystifying,dean2012large}.
Two of the most popular parallel paradigms are data-parallel and model-parallel approaches.
Data-parallel approaches~\cite{krizhevsky2014one} split the training data among multiple processing devices.
Each device then evaluates the loss/gradient for a given portion of data, and the average value is used to update the network parameters.
In contrast, model-parallel approaches aim at splitting the parameters of the DNN model among multiple devices, see for example pipelining techniques~\cite{huang2019gpipe,narayanan2019pipedream} or tensor-parallel approaches~\cite{ben2019demystifying,rudakov2024activations}.

In numerical analysis, domain decomposition (DD) methods are commonly used to accelerate the convergence and parallelize solution strategies for large-scale problems~\cite{dolean2015introduction, toselli2004domain}.
Among these approaches, the additive Schwarz method (ASM)~\cite{lions1988schwarz} is widely utilized for its simplicity and parallel computing capabilities.
In ASM, the computational domain is first decomposed into a collection of smaller, possibly overlapping, subdomains.
The solution of a PDE is then sought in parallel for each subdomain, and the acquired local solutions are subsequently combined to update the global solution approximation.

Although the ASM is parallelizable by design, it is not algorithmically scalable, as its convergence tends to deteriorate with an increased number of subdomains~\cite{dolean2015introduction, toselli2004domain}.
The common remedy for achieving the algorithmic scalability, i.e., keeping the number of iterations bounded independently of the number of subdomains, is to introduce: i) overlap between the subdomains ii) a coarse-level correction step, giving rise to the two-level overlapping additive Schwarz method (TL-ASM)~\cite{dryja1989additive}.
The role of the overlap is to provide the continuity between local solutions at the interface, while the coarse-level correction step ensures the transfer of data across the entire domain at marginal computational cost~\cite{dryja1989towards}.
Building upon these advancements, our goal is to construct the nonlinear TL-ASM in order to parallelize the training of scientific machine learning applications while at the same time enhancing the convergence speed of the current state-of-the-art optimizers.

In the field of machine learning, DD approaches are increasingly utilized to enhance the training of DNNs.
For example, the parallelization of training through the decomposition of deep convolutional networks has been proposed in~\cite{gu2022decomposition,klawonn2023domain}.
In~\cite{gunther2020layer,kopanivcakova2024enhancing,lee2022parareal}, multilevel decompositions giving rise to layer-parallel training have been explored.
Moreover, in the context of PINNs, several DD-based approaches that leverage the spatiotemporal decomposition of the computational domain have been proposed, such as XPINN~\cite{jagtap2021extended}, cPINN~\cite{jagtap2020conservative}, APINN~\cite{hu2022augmented}, and multilevel DD-based architectures~\cite{dolean2023multilevel,LEE2022323}.
In this work, we build upon the recent work presented in~\cite{kopanivcakova2024enhancing}, where single-level additive and multiplicative preconditioners for training of PINNs were introduced.
These preconditioners utilize the layer-wise decomposition of DNNs and are, therefore, agnostic to the underlying physics and computational domain.
This suggests their potential use beyond the training of PINNs to a broader range of scientific machine learning problems, including data-driven operator learning approaches.

Starting from the single-level additive layer-wise preconditioner~\cite{kopanivcakova2024enhancing}, we devise a novel two-level overlapping additive Schwarz preconditioner by leveraging and expanding established techniques from numerical analysis related to nonlinear two-level additive Schwarz methods (TL-ASM)~\cite{cai2002non,dolean2016nonlinear,kopanivcakova2023nonlinear,kothari2022nonlinear_bounds,luo2023preconditioned}.
The main contributions of this work can be summarized as follows.
\begin{itemize}
	\item We introduce overlapping layer-wise decomposition of DNNs, effectively improving the convergence of non-overlapping layer-wise preconditioners.
		\newtext{Motivated by domain-decomposition (DD) approaches in scientific computing, we enforce overlap between the parameters of neighboring layers/subdomains.}
		
	\item  A novel subdomain-wise synchronization strategy is proposed to recombine corrections obtained from parallel subdomain training. This expands upon commonly used methodologies in numerical analysis, which combine subdomain corrections using the weighted averaging technique, see for example literature on the (nonlinear) restricted additive Schwarz method~\cite{cai2002non,chaouqui2022linear}.
	\newtext{In the context of the layer-wise DD methods considered in this work, introducing a subdomain-wise synchronization strategy improves the convergence speed of the overlapping DD method as the size of the overlap increases.}

	\item We incorporate a coarse-level training step that enables the global communication between subdomains while maintaining the feed-forward nature of feature flow in the DNNs. 
	\newtext{This leads to significantly faster convergence compared to using single-level domain decomposition methods.}
	\item  The applicability of the single-level and eventually two-level layer-wise preconditioners is extended beyond PINNs to DeepONets, \newtext{which demonstrate the applicability of the proposed method to wider range of machine-learning problems.}

	\item Through a series of numerical experiments, we demonstrate that the proposed two-level preconditioner yields more accurate DNN models, while also improving the convergence rate of the standard LBFGS optimizer and the baseline single-level preconditioner~\cite{kopanivcakova2024enhancing}.

\end{itemize}

This paper is organized as follows.
\Cref{sec:sciml} provides an introduction to scientific machine learning, with a particular focus on PINNs and DeepONets.
In~\Cref{sec:preconditioner}, we propose a novel two-level overlapping additive Schwarz preconditioner for training DNNs.
\Cref{sec:examples} presents a set of benchmark problems and discusses the implementation details.
Finally, in~\Cref{sec:results}, we demonstrate the numerical performance of the proposed preconditioner.
The conclusion and future work are discussed in~\Cref{sec:conclusion}.

\section{Scientific machine learning}
\label{sec:sciml}
Scientific machine learning is an emerging research field that integrates scientific computing with machine learning techniques to solve complex scientific and engineering problems.
In this work, we focus our attention on solving the parameterized PDEs using PINNs~\cite{raissi2019physics} and DONs~\cite{lu2021learning}.
To this aim, let $\Omega \subset \R^{d}$, $d=1,2$, or $3$, be a closed bounded domain and let the parameter set $\Hv$ be a compact subset of $\R^{p}$, where $p \geq 0$.
The parametric PDE is given in its abstract form as: For $\etav \in  \Hv$, find $u(\xv, \etav) \colon \Omega \times \Hv \to \R$, such that
\begin{equation}
	\label{eq:pde}
	\begin{split}
		\A (u(\xv, \etav), \etav)     & = f(\xv, \etav), \quad \  \xv \in \Omega,                                                                               \\
		\B^{k} (u(\xv, \etav), \etav) & = g^{k}(\xv, \etav), \quad \xv \in \Gamma^{k} \subseteq \partial \Omega, \quad \text{ for } k=1, 2, \ldots, n_{\Gamma},
	\end{split}
\end{equation}
where $\A$ and $\{\B\}_{k=1}^{n_{\Gamma}}$ denote a nonlinear differential operator and a set of boundary condition operators\footnote{In the case of time-dependent problems, we handle the time dimension as an additional component of the vector $\xv \in \Omega$, which allows us to treat the initial condition as a specific type of the boundary conditions.}, respectively.
Moreover, the right hand side function $f(\xv, \etav)$ and the functions $\{g^{k}(\xv, \etav)\}_{k=1}^{n_{\Gamma}}$ are given.

\subsection{Physics informed neural networks}
\label{sec:piml}
PINNs~\cite{raissi2019physics} leverage empirical data and knowledge of physical priors to solve PDEs using DNNs.
Here, we focus on solving one particular instance of~\eqref{eq:pde} with a preselected set of parameters~$\etav$.
A solution $u(\xv, \etav)$ of a PDE for a given $\etav$ at any $\xv \in \Omega$ is approximated by the DNN $u_{\thetav} \colon \R^{d} \to \R$, which is parameterized with respect to parameters $\thetav \in \R^n$.
An actual form of the DNN depends on the user-specified neural network architecture.
For instance, $u_{\thetav}$ can be represented by a feed-forward network or a residual network~\cite{he2016deep}.

In order to ensure that the network~$u_{\thetav}$ approximates the PDE, we have to find the optimal set of DNN's parameters~$\thetav^{\ast}$.
To this aim, we minimize the loss $\L \colon \R^n \to \R$ associated with a residual of a given PDE and boundary conditions, i.e.,
\begin{equation}
	\label{eq:pinn}
	\thetav^{\ast} = \argmin_{\thetav \in \R^n} \L (\thetav) := w_{\mathrm{int}}\L_{\mathrm{int}}(\thetav) + w_{\mathrm{bc}}\L_{\mathrm{bc}}(\thetav),
\end{equation}
where
\begin{equation}
	\label{eq:pinn2}
	\begin{split}
		\L_{\mathrm{int}}(\thetav) & := \Vert \A (u_{\thetav}) - f \Vert^{2}_{L^{2}(\Omega)} = \int_{\Omega} (\A (u_{\thetav}(\xv))-f(\xv))^{2} d \xv,                   \\
		\L_{\mathrm{bc}}(\thetav)  & := \Vert \B (u_{\thetav}) - g \Vert^{2}_{L^{2}(\partial \Omega)} = \int_{\partial \Omega} (\B (u_{\thetav}(\xv))-g(\xv))^{2} d \xv.
	\end{split}
\end{equation}

To discretize the interior loss~$ \L_{\mathrm{int}}$, we consider a set of points $\D_{\mathrm{int}}=\{ \xv_{j}\}^{n_{\mathrm{int}}}_{j=1}$ sampled from the interior of the domain $\Omega$.
Similarly, in order to discretize the boundary loss~$ \L_{\mathrm{bc}}$, we construct $\{\D^{k}_{\mathrm{bc}}\}_{k=1}^{n_{\Gamma}}$, where each point of $\D^{k}_{\mathrm{bc}}=\{ \xv^{k}_{j} \}^{n_{\mathrm{bc}}}_{j=1}$ is sampled from the boundary $\Gamma^k$.
The continuous losses~\eqref{eq:pinn2} can be then discretized, for example, by using the Monte-Carlo method, giving rise to
\begin{equation*}
	\begin{split}
		\L_{\mathrm{int}}(\thetav) & = \frac{1}{\vert \D_{\mathrm{int}}\vert} \sum_{\xv_{j} \in \D_{\mathrm{int}}} (\A (u_{\thetav}(\xv_{j}))-f(\xv_{j}))^{2},                                                          \\
		\L_{\mathrm{bc}}(\thetav)  & = \sum_{k=1}^{n_{\Gamma}} \left( \frac{1}{\vert \D^{k}_{\mathrm{bc}}\vert} \sum_{\xv^{k}_{j} \in \D^{k}_{\mathrm{bc}}} (\B\newtext{^{k}} (u_{\thetav}(\xv^{k}_{j}))-g\newtext{^{k}}(\xv^{k}_{j}))^{2} \right).
	\end{split}
\end{equation*}

The weights $w_{\mathrm{int}}$ and $w_{\mathrm{bc}}$ in~\eqref{eq:pinn} are used to balance the interior and boundary loss terms.
In practice, the quality of the DNN approximation and the success of the training process depends on the appropriate selection of these weights~\cite{sukumar2022exact}.
Consequently, many approaches for manual or automatic selection of weights have been proposed in the literature; see, for instance, NTK-based weighting strategies~\cite{wang2022and} or residual-based attention mechanism~\cite{anagnostopoulos2024residual}.
\newtext{In this work, Dirichlet boundary conditions are directly embedded into the network architecture, thereby avoiding explicit selection of $w_{\mathrm{int}}$ and $w_{\mathrm{bc}}$.}
Following~\cite{sukumar2022exact,lu2021physics}, we modify the network output as
\begin{equation}
	\label{eq:exact}
	\tilde{u}_{\thetav}(\xv) := \sum_{k=1}^{n_{\Gamma}} g^{k}(\xv) + \ell(\xv) u_{\thetav}(\xv), \quad \text{ where } \ell(\xv) :=\prod_{k=1}^{n_{\Gamma}}\ell^{k}(\xv).
\end{equation}
Each function $\ell^{k} \colon \Omega \to [0,1] $ is defined such that
\begin{equation*}
	\ell^{k}(\xv) = \begin{cases}
		0,           & \xv \in \Gamma^{k},                                                                      \\
		1,           & \xv \in \Gamma^{j}, \quad j \in \{ 1,2,\ldots,n_{\Gamma} \}, \quad \text{ for } k \ne j, \\
		t \in (0,1), & \text{otherwise}.
	\end{cases}
\end{equation*}
Note that $t$ is determined so that $\ell^{k}(\xv)$ is a smooth function.
For example, when $\Omega = (0,1)^{2}$ and $\partial \Omega = \cup_{k=1}^{n_{\Gamma}} \Gamma^{k}$, $\ell(\xv)=\ell(x_{1}, x_{2})=x_{1}(1-x_{1})x_{2}(1-x_{2})$.
Using~$\tilde{u}_{\thetav}$ as defined in~\eqref{eq:exact}, we can now reformulate the minimization problem~\eqref{eq:pinn} as
\begin{equation}
	\thetav^{\ast} = \argmin_{\thetav \in \R^n} \L (\thetav) := \frac{1}{\vert \D_{\mathrm{int}} \vert} \sum_{\xv_{j} \in \D_{\mathrm{int}}} (\A (\tilde{u}_{\thetav}(\xv_{j})) - f(\xv_{j}) )^{2}.
	\label{eq:loss_pin_exact_bc}
\end{equation}

\subsection{Operator learning}
\label{sec:ol}
Operator learning approaches use DNNs to approximate a mapping between infinite-dimensional function spaces.
Let $\Y$ be an infinite-dimensional Banach space, which can represent a space of parametrized right-hand sides, boundary conditions, or material parameters.
The goal is to approximate the operator $\G \colon \Y \to \V$, where $\V$ denotes the solution space of the parametric PDE, specified by~\eqref{eq:pde}.
In this work, we approximate the operator $\G$ using the DON, initially proposed in~\cite{lu2021learning}.

The architecture of DON consists of branch and trunk networks.
The branch network $B \colon \Y_{m} \to \R^{p}$ provides a set of  $p$ coefficients for a given discretized function ${\yv\newtext{^{m}}\oldtext{_{m}} \in \Y_{m}}$.
Here,   $\Y_{m}$ is a $m$-dimensional space that is assumed to satisfy the following condition:
\begin{equation*}
	\forall y \in \Y, \  \exists \yv\newtext{^{m}}\oldtext{_{m}} \in \Y_{m}\newtext{\subseteq \Y}  \text{ such that } \yv\newtext{^{m}}\oldtext{_{m}} \to y \text{ as } m \to \infty.
\end{equation*}
\newtext{For simplicity, we write $\yv$ in place of $\yv^{m}$, as $m$ is a fixed and determines the dimension of the discrete space.}
In addition, the trunk network~$T \colon \Omega \to \R^{p}$ is used to encode the computational domain~$\Omega$.
Therefore, it takes as an input a coordinate $\xiv \in \Omega$ and returns a set of $p$ basis functions.

\begin{figure}
	\centering
	\includegraphics[scale=1.1]{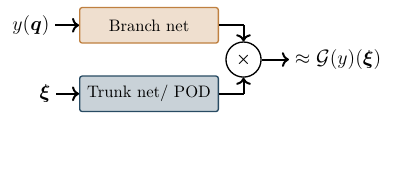}
	\hfill
	\includegraphics[scale=1.0]{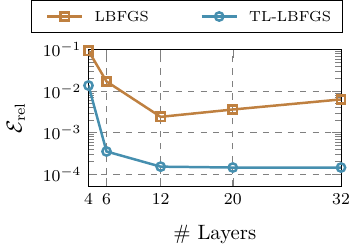}
	\caption{Left: A sketch of DON architecture.
		The trunk can be represented by the DNN or by the POD basis obtained from the training data.
		The symbol $\boldsymbol{\times}$ denotes the inner product between the output of the branch and trunk networks.
		Right: A relative error obtained by using the standard LBFGS optimizer and the proposed TL-LBFGS as a function of increasing number of layers. The experiment is performed for the AC example, specified in~\Cref{sec:PINNS_examples}.}
	\label{fig:deeponet}
\end{figure}

The coefficients $B_{k}(\yv) \in \R^p$ and the basis functions $T_{k}(\xiv) \in \R^p$ obtained as an output of branch and trunk networks are then combined as
\begin{equation}
	\G_{\thetav} (\yv) (\xiv) := \sum_{k=1}^{p} B_{k}(\yv) T_{k}(\xiv), \quad \yv \in \Y_{m}, \xiv \in \Omega,
\end{equation}
to approximate a solution of the parametric PDE at a point~$\xiv$ for a given function~$\yv$.
Figure~\ref{fig:deeponet} illustrates a prototypical DON architecture.
Note that we are free to choose the architecture of branch and trunk networks.
For example, the choice of the branch network can be guided by the structure of the input~$\yv$, e.g., it can be a fully connected network in case of a scalar input, or a convolutional network in case of multidimensional discretized functions.
Since the vector of coordinates $\xiv \in \R^d$ is usually low dimensional, fully connected networks are commonly used in place of the trunk network.
Alternatively, one can utilize the precomputed basis functions, e.g., obtained by means of the proper orthogonal decomposition~(POD) performed on the training data~\cite{lu2022comprehensive} or by means of parameterizing the integral kernel in the Fourier space~\cite{anandkumar2019neural}.

To train DON such that it approximates solutions of~\eqref{eq:pde}, we construct a dataset
\begin{equation*}
	\D = \{ (\yv_{j}, \bar{\xiv}_{j}, \uv_{j}) \}_{j=1}^{n_{s}},
\end{equation*}
consisting of $n_{s}$ samples.
Each triple~$(\yv_{j}, \bar{\xiv}_{j}, \uv_{j})$ consists of the discretized function ${\yv_{j} \in \R^m}$,  sampled from $\Y_{m}$.
Moreover,  $\bar{\xiv}_{j} = [ \xiv_{j,1}, \ldots, \xiv_{j, n_{c}} ]^T \in \R^{n_{c} \times d}$ represents a set of $n_{c}$ collocation points and $\uv_{j} \in \R^{n_{c}}$ is a target solution evaluated at collocation points constituting $\bar{\xiv}_{j}$.
\newtext{Note that the target solution is typically obtained using high-fidelity discretization techniques, such as finite element or finite difference methods.}

Using dataset~$\D$, the optimal parameters of the DON, denoted by $    \thetav^{\ast} \in \R^n$, are found by solving the following minimization problem:
\begin{equation}
	\label{eq:donloss}
	\thetav^{\ast} = \argmin_{\thetav \in \R^n}    \L (\thetav) := \frac{1}{n_{s} n_{c}} \sum_{j=1}^{n_{s}} \sum_{t=1}^{n_{c}} \left(\G_{\thetav}(\yv_{j})(\xiv_{j,t}) - \uv_{j}\right)^{2}.
\end{equation}
Note that instead of using $\uv_{j}$ given by the dataset, we can utilize the PINN-based loss function; see~\cite{wang2021learning} for details.
This might simplify the generation of the dataset but makes the training process more tedious.
We do not pursue this possibility here but rather investigate the performance of the proposed two-level preconditioners in data-driven settings.
Notably, the proposed preconditioner is expected to enhance the performance of standard optimizers more significantly in the case of PINN-based loss due to increased ill-conditioning of the training problem~\cite{wang2021learning}.

\section{Two-level overlapping additive Schwarz preconditioner}
\label{sec:preconditioner}
In this work, we focus on training the DNN by solving the following minimization problem:
\begin{align}
	\min_{\thetav} \L (\thetav),
	\label{eq:min_problem}
\end{align}
where the loss $\L$ is as defined in~\eqref{eq:loss_pin_exact_bc}  or~\eqref{eq:donloss} for PINN or DON, respectively.
The minimization problem~\eqref{eq:min_problem} is traditionally solved using gradient-based optimizers such as Adam~\cite{kingma2014adam} or LBFGS~\cite{liu1989limited}.
However, the convergence of these optimizers tends to deteriorate as problem stiffness and ill-conditioning increase~\cite{wang2021understanding}, resulting in slow training and unsatisfactory accuracy of the resulting PINN/DON models.
For instance,~\Cref{fig:deeponet} illustrates that increasing the number of layers in a PINN model exacerbates \newtext{the relative error when using a standard optimizer.
As the depth of the network grows, the ill-conditioning and non-convexity of the underlying minimization problem become more significant.}
This, in turn, presents significant challenges for the standard LBFGS method to converge, leading to poorer accuracy of the PINN as the representation capacity of the DNN increases.

\newtext{In scientific computing, Schwarz methods~\cite{lions1988schwarz, matsokin1985schwarz} provide an efficient framework for solving large-scale problems by decomposing the solution space into smaller subspaces, within which the problem can be solved more effectively.}
Emulating the classical additive Schwarz methods\oldtext{To overcome these difficulties}, the single-level Schwarz preconditioning strategy, which utilizes the layer-wise decomposition of the DNN, was introduced in~\cite{kopanivcakova2024enhancing}.
In this work, we propose enhancing the performance of this strategy by introducing overlap, a coarse-level correction step, and a subdomain-based recombination strategy.
The culmination of these techniques results in a novel two-level overlapping Schwarz preconditioner, the performance of which does not deteriorate with increasing network depth, as shown in~\Cref{fig:deeponet}.

To present the proposed two-level overlapping Schwarz preconditioning strategy, let us consider the nonlinear systems of equations associated with the first-order optimality conditions of~\eqref{eq:min_problem} given as
\begin{equation}
	\label{eq:fcond}
	\nabla \L (\thetav) = 0.
\end{equation}
Since solving~\eqref{eq:fcond}  is computationally demanding,  we instead construct and solve the following nonlinearly-preconditioned system of equations:
\begin{equation}
	\label{eq:precond}
	\F (\thetav) := \nabla \L (\pazocal{M}(\thetav)) = 0,
\end{equation}
where the preconditioner $\pazocal{M} \colon \R^{n} \to \R^{n}$ is designed such that~\eqref{eq:precond} has more balanced nonlinearities than~\eqref{eq:fcond}. 
\newtext{In practice, this means that the nonlinear terms are of comparable magnitude both to each other and to the linear terms, thereby avoiding dominance by any single contribution.}

In this work, we solve the nonlinear system~\eqref{eq:precond} iteratively.
Thus, on each $k$-th iteration, we start from a given iterate $\thetav^{(k)}$, and apply the preconditioner $\pazocal{M}$.
As a result, we obtain the \newtext{intermediate}\oldtext{smoothed} iterate $\thetav^{(k+1/2)} := \pazocal{M}(\thetav^{(k)})$, which is required to reduce the loss, i.e., it has to satisfy the following condition:
\begin{equation*}
	\L (\thetav^{(k+1/2)}) \leq \L (\thetav^{(k)}).
\end{equation*}
The smoothing step is followed by a global optimization step.
To this aim, we apply an operator~$\K \colon \R^n \to \R^n$, which is associated with one or multiple iterations of a user-specified global optimizer.
The operator~$\K$  takes as an input the gradient of the loss function and provides a direction $\pv^{(k)}$, i.e.,~${\pv^{(k)} := \K (\nabla \L (\thetav^{(k+1/2)}))}$.
For instance, if a single step of the LBFGS optimizer is employed, then the operator $\K$ is defined as
\begin{equation*}
	\pv^{(k)} :=     \K \big(\nabla \L ( \thetav^{(k+1/2)}) \big) := - \bigg( \Hm^{(k+1/2)} \bigg)^{-1} \nabla \L(\thetav^{(k+1/2)}),
\end{equation*}
where $\Hm^{(k+1/2)}$ denotes the Hessian approximation of $\L$ at $\thetav^{(k+1/2)}$.
Afterward, the updated global iterate $\thetav^{(k+1)}$ is obtained as
\begin{equation*}
	\thetav^{(k+1)} = \thetav^{(k+1/2)} + \alpha^{(k)}\pv^{(k)},
\end{equation*}
where $\alpha^{(k)}$ denotes the step size, determined using the line-search globalization strategy.

\subsection{Overlapping layer-wise additive Schwarz preconditioner}
\label{sec:net_decomposition}
\begin{figure}
    \includegraphics{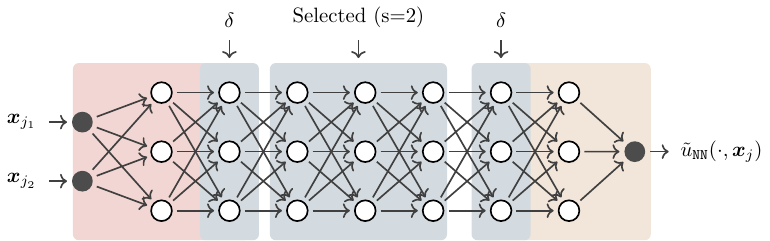}
	\caption{An illustration of the network decomposed into three subdomains with overlap $\delta=1$. 
    The parameters of different color belong to different subdomains. 
    The second subdomain (blue color, $s=2$) is selected to illustrate the overlapping subdomain parameters.}
	\label{fig:overlap_decomp_net}
\end{figure}

The computational \newtext{efficiency}\oldtext{cost} associated with solving the preconditioned system of equations~\eqref{eq:precond} is directly linked to the quality \newtext{and the computational cost }of the nonlinear preconditioner~$\pazocal{M}$.
In this work, we construct the preconditioner~$\pazocal{M}$ by utilizing the overlapping layer-wise network decomposition.
Thus, we decompose the DNN, consisting of $N$ layers, into $N_{sd}$ subnetworks (subdomains), such that each layer belongs to a unique subnetwork.
The parameters of the network are also decomposed into $N_{sd}$ disjoint groups, i.e.,
\begin{equation*}
	\thetav = [\thetav_{1}, \ldots, \thetav_{s}, \ldots, \thetav_{N_{sd}}]^{T},
\end{equation*}
where $\thetav_{s} := [\thetav_{s, 1}, \ldots, \thetav_{s, N_s}]^T \in \R^{n_{s}}$, for $s = 1, \ldots, N_{sd}$.
Each $s^{th}$ non-overlapping subnetwork consists of $N_s$ layers, with parameters of $h^{th}$ layer being denoted as~$\thetav_{s, h}$.
In the case of the overlapping decomposition, the parameters $\thetav_{s}$ of the $s$-th subnetwork are concatenated with the parameters of $\delta(\in \N)<N_S$ neighboring layers from the preceding $(s-1)$-th and succeeding $(s+1)$-th subnetworks.
For~$\delta\geq 1$, this gives rise to overlapping set of subdomain parameters
\begin{align*}
	\tilde{\thetav}_{s} = [{\thetav}_{s-1, (N_{sd}-\delta)}, \ldots,  {\thetav}_{s-1, N_{sd}},  {\thetav}_{s, 1}, \ldots, {\thetav}_{s, N_{sd}}, {\thetav}_{s+1, 1}, \ldots,  {\thetav}_{s+1,  \delta}]^{T} \in \R^{\tilde{n}_{s}},
\end{align*}
where $\tilde{N}_{s}$ and $\tilde{n}_{s}$ denotes the number of layers in overlapping subnetwork and the number of associated subnetwork parameters, respectivelly.
Note that, for the first ($s=1$) and the last ($s = N_{sd}$-th) subnetworks, we only pick layers from the second and $(N_{sd}-1)$-th subnetworks, respectively.
\Cref{fig:overlap_decomp_net} illustrates the parameter selection for the second subnetwork of the six-layer network decomposed in a layer-wise manner into three subdomains with an overlap of size one, i.e.,~$N_{sd}=3$ and $\delta=1$.

\newtext{
\begin{remark}
The subdomain size can range from a minimum of $1$ (serial) to a maximum of $n$ layers.  
The overlap size can range from a minimum of $0$ (non-overlapping) to the maximum number of layers in a subnetwork ($n_{s} = n / N_{sd}$ if the layers are decomposed uniformly).  
\end{remark}}

To transfer the data between the subnetworks and the global network, we define a restriction operator $\Rv_{s} \colon \R^{n} \to \R^{\tilde{n}_{s}}$, which extracts the parameters of the overlapping subnetwork from the parameters of the global network, i.e.,
\begin{equation*}
	\tilde{\thetav}_{s} = \Rv_{s} \thetav.
\end{equation*}
Moreover, we also consider the prolongation operator $\Rv^{T}_{s} \colon \R^{\tilde{n}_{s}} \to \R^{n}$, defined as adjoint of $\Rv_{s}$.
The role of the prolongation operator $\Rv^{T}_{s}$ is to extend the parameters associated with $s^{th}$ subnetwork to the global network, i.e., $\thetav = \sum_{s=1}^{N_{sd}} \Rv^T_{s}  \tilde{\thetav}_{s}$.

\begin{remark}
	In the case of DON, we decompose the trunk and branch networks independently of each other.
	For instance, if $N_{sd}=2$, the branch and trunk constitute two separate non-overlapping subdomains.
	If $N_{sd}>2$, then branch and trunk networks are further decomposed into multiple, possibly overlapping, subdomains.
\end{remark}

\subsubsection{Single-level preconditoner with layer-wise synchronization strategy}
\label{sec:recomination_strategy}
Using the proposed overlapping network decomposition, we define the local training problems for all subdomains as
\begin{equation}
	\label{eq:local}
	\tilde{\thetav}^{\ast}_{s} := \argmin_{\tilde{\thetav}_{s} \in \R^{n_s}} \L (\tilde{\thetav}_{s} + (\Im - \Rv^{T}_{s}\Rv_{s})\thetav),
\end{equation}
where $s \in 1, \ldots, N_s$.
Thus, for the $s^{th}$ subdomain,  we minimize the loss functional $\L$ with respect to parameters~$\tilde{\thetav}_{s}$,  while all other parameters are held fixed.
The minimization problem~\eqref{eq:local} is solved only approximately, for example, by utilizing a few steps of the LBFGS or Adam optimizer.
The resulting parameters~$\tilde{\thetav}^{\ast}_{s}$ are then used to define the (single-level) overlapping additive Schwarz preconditioner $\pazocal{M}_{SL}$ as
\begin{equation}
	\label{eq:onelevel}
	\pazocal{M}_{SL}(\thetav) := \thetav + \gamma \sum_{s=1}^{N_{sd}} \Rv^{T}_{s}(\tilde{\thetav}^{\ast}_{s} - \Rv_{s} \thetav) =: \thetav +  \gamma \sum_{s=1}^{N_{sd}}  \dv_{s} =: \thetav + \gamma \dv.
\end{equation}
Here, the symbol $\gamma \in [0, 1]$ denotes an appropriately selected damping parameter, such that
\begin{equation*}
	\L (\thetav + \gamma \dv) \leq \L (\thetav).
\end{equation*}

It is crucial to note that the convergence rate of the preconditioner $\pazocal{M}_{SL}$ is strongly affected by the value of $\gamma$ as well as by the quality of the search direction $\dv$.
In the literature on the linear additively composed preconditioners, the multipreconditioning approach~\cite{bridson2006multipreconditioned,greif2017gmres} is often employed to determine how to optimally combine locally obtained search directions~$\{\dv_{s}\}_{s=1}^{N_{sd}}$ to improve the quality of the resulting search direction~$\dv$.
Here, we extend this methodology into nonlinear settings by searching for a set of  damping parameters $\{\gamma^*_{s}\}_{s=1}^{N_{sd}}$, such that
\begin{align}
	\L \left(\thetav + \sum_{s=1}^{N_{sd}}\gamma^*_{s} \dv_{s} \right) < \L \left(\thetav + \sum_{s=1}^{N_{sd}} \gamma_{s} \dv_{s} \right),\quad \newtext{ \forall \gamma^*_s \neq \gamma_s,}
\end{align}
\newtext{where $\gamma_{s} \in [0, 1]$ is an initial damping parameter associated with the local, subdomain-based, search direction $\dv_{s}$.}\oldtext{for all~$\gamma^*_s \neq \gamma_s$. 
In this work, we choose initial damping parameter as $\gamma_{s}=1$.}

The optimal set of parameters $\{\gamma^*_{s}\}_{s=1}^{N_{sd}}$ can be determined by solving the following $N_{sd}$-dimensional minimization problem:
\begin{align}
	\min_{\gamma_{1}, \ldots, \gamma_{N_{sd}}}  \L\left(\thetav +  \sum_{s=1}^{N_{sd}}\gamma_{s} \dv_{s} \right).
	\label{eq:gammas_min_problem}
\end{align}
However, finding a solution of~\eqref{eq:gammas_min_problem} is computationally demanding due to the nonlinear nature of the problem.
This is in contrast to multi-preconditioning approaches used in the context of linear problems, where $\{\gamma^*_{s}\}_{s=1}^{N_{sd}}$ can  typically be obtained in closed form, see for instance~\cite{bridson2006multipreconditioned,brezinski1999multiparameter}.

To elevate the computational burden of solving~\eqref{eq:gammas_min_problem} exactly, one can employ multidimensional line-search heuristics, such as multidimensional backtracking proposed in~\cite{kunstner2024searching}.
Here, we opt for a simplified approach and leverage the fact that the majority of DNNs construct their bases by propagating the input features through the network from the first to the last layer in a sequential manner, which is frequently regarded as a time-dependent process~\cite{gunther2020layer,lee2022parareal}.
Therefore, we also follow the feature flow and determine $\{\gamma_{s}^*\}_{s=1}^{N_{sd}}$ in a sequential manner, such that
\begin{equation*}
	\L \left( \thetav + \sum_{j=1}^{s-1} \gamma_{j}^* \dv_{j}  + \gamma_{s} \dv_{s}  \right) \leq \L \left( \thetav + \sum_{j=1}^{s-1} \gamma_{j}^* \dv_{j} \right),
\end{equation*}
for $s=1, \ldots, N_{sd}$.
Thus, the values~$\{\gamma^*_{s}\}_{s=1}^{N_{sd}}$ are determined from the first subdomain to the last one, utilizing a standard, unidimensional, line-search strategy.
This guarantees that each search direction~$\dv_{s}$ decreases the value of $\L$ associated with a global network.
At the same time, we ensure that the DNN's bases are constructed by appropriately adjusting the impact of each search direction~$\dv_{s}$ on the global DNN's approximation.
Here, we remark that the search directions $\{\dv_{s}\}_{s=1}^{N_{sd}}$ are obtained in parallel, while only scalar damping parameters $\{\gamma^*_{s}\}_{s=1}^{N_{sd}}$ are obtained sequentially during the subdomain synchronization step.

\subsection{Two-level overlapping additive Schwarz preconditioner}
\label{sec:coarse_level}
\begin{figure}
 \includegraphics{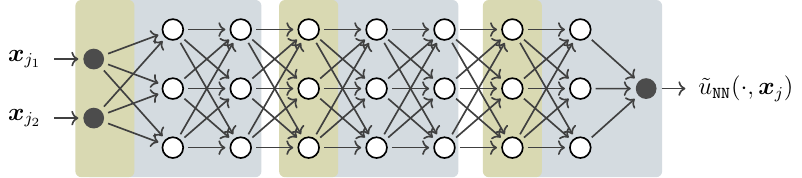}
	\caption{An example of network decomposition associated with coarse level training.
	The first layer (olive color) of each subdomain (olive+blue color) is selected to construct a set of trainable coarse level parameters.}
	\label{fig:coarse}
\end{figure}

In the field of numerical methods, the convergence of single-level preconditioners is often enhanced by incorporating a coarse level, which ensures a global communication between subdomains.
Following this methodology, we propose to improve the convergence of a single-level layer-wise preconditioner by incorporating a coarse-level training step.
Our design of the coarse-level network is inspired by the observation that the forward pass through the network can be interpreted as a time-stepping process~\cite{gunther2020layer,lee2022parareal,gaedke2021multilevel,gratton2023multilevel,kopanicakova2022globally}.
This observation suggests that coarse-level subnetworks can be created by reducing the number of layers within a given subdomain, analogous to time-step coarsening.
More precisely,  we create the coarse-level subnetwork by grouping together the first layers of each subdomain.
Thus, the trainable coarse-level parameters are given as
\begin{align}
	\thetav_{0} = [\thetav_{1, 1}, \ldots, \thetav_{s, 1}, \ldots, \thetav_{N_{sd}, 1} ]^T,
\end{align}
where~$\thetav_{s, 1}$ is extracted from $\thetav_{s} = [\thetav_{s, 1}, \ldots, \thetav_{s, \tilde{N}_s}]^T$, for all $s=1, \ldots, N_{sd}$.
An illustration of the coarse-level parameter selection for a six-layer network with  $N_{sd}=3$ and $\delta=0$ is depicted in~\Cref{fig:coarse}.
\newtext{The computational cost of the coarse-level problem is proportional to the number of subdomains $N_{sd}$. }

The transfer of the data between the coarse-level subnetwork and the global network, is ensured by two transfer operators.
In particular, we consider the restriction operator $\Rv_{0} \colon \R^{n} \to \R^{n_{0}}$,\oldtext{ and its adjoint, the   prolongation operator~$\Rv^{T}_{0} \colon \R^{n_{0}} \to \R^{n}$} given such that
\begin{equation*}
	\thetav_{0} = \Rv_{0} \thetav\oldtext{, \qquad \qquad \thetav = \Rv_{0}^T \thetav_{0}}.
\end{equation*}
\newtext{The prolongation operator~$\Rv^{T}_{0} \colon \R^{n_{0}} \to \R^{n}$ is defined as the adjoint operator of $\Rv_{0}$.}
Note that the number of trainable coarse-level parameters~$n_{0}$ scales proportionally to the number of subdomains $N_{sd}$.

Using the coarse-level subnetwork, we can now define the coarse-level training problem as
\begin{equation*}
	\thetav^{\ast}_{0} := \argmin_{\thetav_{0} \in \R^{n_0}} \L (\thetav_{0} + (\Im - \Rv^{T}_{0}\Rv_{0})\thetav).
\end{equation*}
Hence, the loss functional $\L$ is minimized with respect to the parameters associated with the coarse-level subnetwork while all other parameters are held fixed.
The outcome of the coarse-level minimization process, the parameters $\thetav^{\ast}_{0} \in \R^{n_0}$, is then used to update the parameters of the global network as
\begin{equation}
	\label{eq:coarse}
	\hat{\thetav} = \thetav + \Rv^{T}_{0} (\thetav^{\ast}_{0} - \Rv_{0}\thetav).
\end{equation}
Note that since all coarse-level parameters $\thetav_{0}$ are optimized simultaneously, the coarse-level training process is serial. Moreover, the condition $\L (\hat{\thetav}) \leq \L (\thetav)$ is automatically fulfilled.

\begin{figure}
\includegraphics{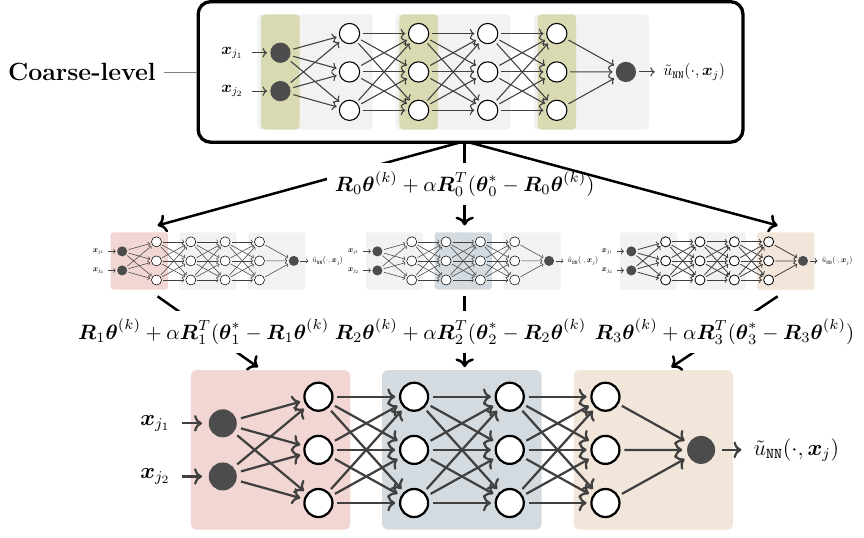}
\caption{\newtext{Schematic illustration of the two-level DD optimizer for training a DNN decomposed into three subdomains.  
We first train the coarse-level parameters (shown in olive color), followed by single-level parallel training of the subdomains.  
The trainable parameters of each subdomain are depicted in different colors, while the non-trainable (fixed) parameters are shown in gray.
In the end, the global training step is performed by performing a few iterations of the optimizer $\pazocal{K}$ for all parameters of the DNN.}}
	\label{fig:schwarz}
\end{figure}

Using~\eqref{eq:coarse}, we are now ready to define a two-level overlapping additive Schwarz preconditioner $\pazocal{M}_{TL}: \R^{n} \to \R^n$, by multiplicatively composing the coarse-level training step with an application of  $\pazocal{M}_{SL}$.
In particular, $\pazocal{M}_{TL}$ is defined by applying $\pazocal{M}_{SL}$ after the coarse-level update~\eqref{eq:coarse} has been performed, i.e.,
\begin{equation*}
	\pazocal{M}_{TL} (\thetav) := \pazocal{M}_{SL} ( \underbrace{\thetav + \Rv^{T}_{0} (\thetav^{\ast}_{0} - \Rv_{0}\thetav)}_{:=\hat{\thetav}}).
\end{equation*}
Note that since the coarse-level training indirectly enforces the feed-forward flow of the features through the network, the parameters~$\hat{\thetav}$ serve as an improved initial guess for the parallel subdomain training phase.
\newtext{\Cref{fig:schwarz} provides a graphical illustration of the one-level and two-level Schwarz methods when applied to train the DNNs.}
In addition, \cref{alg:tlspqn} summarizes the proposed two-level preconditioning process for a user-selected optimizer~$\pazocal{K}$.
\newtext{Our practical implementation of the line-search step (Step~9 in \cref{alg:tlspqn}) employs cubic backtracking with the Strong-Wolfe conditions~\cite{nocedal2006numerical}.}

\begin{algorithm}
	\caption{Two-level additive Schwarz preconditioner for $\K$ optimizer (TL-$\K$)}\label{alg:tlspqn}
	\begin{algorithmic}[1]
		\REQUIRE $\L \colon \R^{n} \to \R$, initial parameters $\thetav^{(0)} \in \R^{n}$, $\varepsilon \in \R$, $k=0$, $\delta \in \N$, $\K:\R^{n} \to \R^n$
		\WHILE{$\L (\thetav^{(k)}) > \varepsilon$}
		\STATE $\thetav^{\ast}_{0} = \argmin_{\thetav_{0}} \L (\thetav_{0} + (\Im - \Rv^{T}_{0}\Rv_{0})\thetav^{(k)})$ \hfill\COMMENT{Coarse level training}
		\STATE $\hat{\thetav}^{(k)} = \thetav^{(k)} + \Rv^{T}_{0} (\thetav^{\ast}_{0} - \Rv_{0}\thetav^{(k)})$
		\FOR{$s=1,\ldots,N_{sd}$}
		\STATE $\tilde{\thetav}^{\ast}_{s} = \argmin_{\tilde{\thetav}_{s}} \L (\tilde{\thetav}_{s} + (\Im - \Rv^{T}_{s}\Rv_{s})\hat{\thetav}^{(k)})$ \hfill\COMMENT{Parallel subdomain training}
		\ENDFOR
		\FOR{$s=1, \ldots, N_{sd}$}
		\STATE $\dv_{s} = \Rv^{T}_{s}(\tilde{\thetav}^{\ast}_{s} - \Rv_{s} \hat{\thetav}^{(k)})$
		\STATE Find $\gamma_{s}$ that satisfies $\displaystyle \L \left( \hat{\thetav}^{(k)} + \sum_{j=1}^{s} \gamma_{j} \dv_{j} \right) \leq \L \left( \hat{\thetav}^{(k)} + \sum_{j=1}^{s-1} \gamma_{j} \dv_{j} \right)$
		\ENDFOR
		\STATE $\displaystyle \thetav^{(k+1/2)} = \hat{\thetav}^{(k)} + \sum_{s=1}^{N_{sd}} \gamma_{s} \dv_{s}$ \hfill \COMMENT{Subdomain correction recombination step}
		\STATE $\pv^{(k)} = \K \big( \nabla \L (\thetav^{(k+1/2)}) \big)$ \hfill\COMMENT{Global optimizer step}
		\STATE $\thetav^{(k+1)}=\thetav^{(k+1/2)}+\alpha^{(k)}\pv^{(k)}$ \hfill\COMMENT{Global parameter update}
		\STATE $k \gets k + 1$
		\ENDWHILE
	\end{algorithmic}
\end{algorithm}

\section{Benchmark problems and implementation details}
\label{sec:examples}
This section describes the benchmark problems used to evaluate the performance of the proposed two-level preconditioner.
In particular, we consider three benchmark problems associated with training PINNs and three benchmark problems with training DONs.
Additionally, we discuss the implementation aspects and configuration of all considered optimizers.

\subsection{Benchmark problems for training of PINNs}
\label{sec:PINNS_examples}
For training PINNs, we consider the following benchmark problems.
\begin{itemize}[leftmargin=0.5cm]
	\item{\textbf{Burgers' equation (Burg):}}
	      Burgers' equation is given in the following form:
	      \begin{equation*}
		      \begin{aligned}
			      \frac{\partial u}{\partial t} + u u' - \nu u'' & = 0,  \quad \quad \quad &  & \forall \ (t,x) \in (0,1] \times (-1,1), \\
			      u                                              & = - \sin(\pi x), \      &  & \forall \ (t,x) \in \{0\} \times [-1,1], \\
			      u                                              & = 0, \                  &  & \forall \ (t,x) \in (0,1] \times \{1\},  \\
			      u                                              & = 0, \                  &  & \forall \ (t,x) \in (0,1] \times \{-1\},
		      \end{aligned}
	      \end{equation*}
	      where~$u=u(t,x)$ denotes the flow velocity and~$\nu$ stands for the kinematic viscosity.
	      Here, we choose~$\nu=0.01/\pi$.

	\item{\textbf{Diffusion-advection (DA):}}
	      The diffusion-advection equation is given as
	      \begin{equation*}
		      \begin{aligned}
			      -\mu \Delta u  + \boldsymbol{b}\cdot \nabla u & = f, &  & \forall \ (x_1,x_2) \in (0, 1) \times (0, 1), \\
			      u                                             & = 0, &  & \text{on } \partial \Omega,
		      \end{aligned}
	      \end{equation*}
	      where~$\boldsymbol{b}=(1,1)^\top$.
	      The right-hand side is considered to be constant, i.e.,~$f = 1$.
	      Moreover, the symbol~$\mu$ denotes viscosity, which we choose to set as~$\mu=10^{-2}$.

	\item{\textbf{Allen-Cahn (AC):}}
	      We consider the Allen-Cahn equation of the following form:
	      \begin{equation*}
		      \begin{aligned}
			      \frac{\partial u}{\partial t}  - D u'' - 5 (u - u^3) & = 0,  \quad \quad \quad &  & \forall \ (t,x) \in (0,1] \times (-1,1), \\
			      u                                                    & = x^2 \cos(\pi x), \    &  & \forall \ (t,x) \in \{0\} \times [-1,1], \\
			      u                                                    & = - 1, \                &  & \forall \ (t,x) \in (0,1] \times \{-1\}, \\
			      u                                                    & = -1, \                 &  & \forall \ (t,x) \in (0,1] \times \{1\},
		      \end{aligned}
	      \end{equation*}
	      where the diffusion coefficient~$D$ is chosen as $D=0.001$.
\end{itemize}

The amount and location of collocation points significantly influence the discretization error of PINNs.
Here, we perform sampling by employing the Quasi-Monte Carlo method~\cite{morokoff1995quasi}, with  Hammersley low-discrepancy sequences~\cite{mishra2022estimates}.
As a consequence, it is well known that the convergence rate of Quasi-Monte Carlo is close to $O(N^{-1})$, while the convergence rate of the standard Monte Carlo method, with pseudo-random sequences, is $O(N^{-1/2})$~\cite{asmussen2007stochastic}.
The number of collocation points is $10{,}000$ for all three examples.

\begin{table}[t]
	\centering
	\caption{The summary of network architectures for PINN's benchmark problems.\newtext{ We use the fully connected neural network~(FNN) with hyperbolic tangent~(Tanh) activation function for Burg, DA, and AC problems. The symbol $[ \cdot ]$ indicates the width for FNN.}}
	\label{tab:pinn_architectures}
	\begin{tabular}{|l||r|r|}
		\hline
		{Example} & \multicolumn{1}{c|}{{Layers}}                     & {Act.} \\ \hline \hline
		Burg      & FNN[2, 20, 20, 20, 20, 20, 20, 1]                 & Tanh   \\  \hline
		DA        & FNN[2, 50, 50, 50, 50, 50, 50, 50, 50, 1]         & Tanh   \\  \hline
		AC        & FNN[2, 32, 32, 32, 32, 32, 32, 32, 32, 32, 32, 1] & Tanh   \\  \hline
	\end{tabular}
\end{table}

\subsection{Benchmark problems for training of DONs}
For the training of DONs, we consider the three benchmark problems, with datasets created using the finite element method within the Firedrake library~\cite{FiredrakeUserManual}.

\begin{itemize}[leftmargin=0.5cm]
	\item{\textbf{Anisotropic Poisson equation (Aniso):}}
	      Let $\Omega := (-1,1)^2$  be the computational domain, with boundary $\partial \Omega$, the anisotropic Poisson equation is given as
	      \begin{equation}
		      \begin{aligned}
			      - \nabla \cdot (K(x, \boldsymbol{\eta}) \ \nabla u(x) ) & = f(x),  \quad \quad &  & \forall \ x \in \Omega,           \\
			      u                                                       & = 0, \               &  & \forall \ x \in  \partial \Omega, \\
		      \end{aligned}
		      \label{eq:aniso_laplace}
	      \end{equation}
	      where~$u$ denotes the solution and $f(x):=1$ is the forcing term.
	      The diffusion coefficient~$K(x, \boldsymbol{\eta})$ is considered to be an anisotropic tensor of the following form:
	      \begin{align*}
		      K(x, \boldsymbol{\eta}) =
		      \begin{pmatrix}
			      \cos(\alpha) & -\sin(\alpha) \\
			      \sin(\alpha) & \cos(\alpha)
		      \end{pmatrix}
		      \begin{pmatrix}
			      1 & 0     \\
			      0 & \beta
		      \end{pmatrix}
		      \begin{pmatrix}
			      \cos(\alpha)  & \sin(\alpha) \\
			      -\sin(\alpha) & \cos(\alpha)
		      \end{pmatrix},
	      \end{align*}
	      where~$\boldsymbol{\eta}:=[\alpha, \beta]$.
	      \newtext{The parameter~${\beta \in [10^{-6},1]}$ specifies the anisotropic strength, and it is sampled from the following distribution: ${\log_{10} (\sfrac{1}{\beta}) \sim \pazocal{U}[0, 6]}$.
	      The parameter~${\alpha \in (0,\pi)}$, sampled as~$\alpha \sim \pazocal{U}[0, \pi]$, is used to specify the angle of anisotropic direction.}
          \oldtext{The components of the forcing term are selected randomly from the Gaussian distribution $\pazocal{N}(0,1)$.}
	      To construct the dataset, for each $\boldsymbol{\alpha}$, we discretize~\eqref{eq:aniso_laplace} using a uniform mesh with $34 \times 34$ elements.
	      The number of training samples and test samples are $4{,}250$ and $750$, respectively.
	      We have uploaded this new dataset to the public ZENODO repository~\cite{zenodo_aniso}.

	      \begin{figure}
		      \hfill
		      \includegraphics[scale=0.0325]{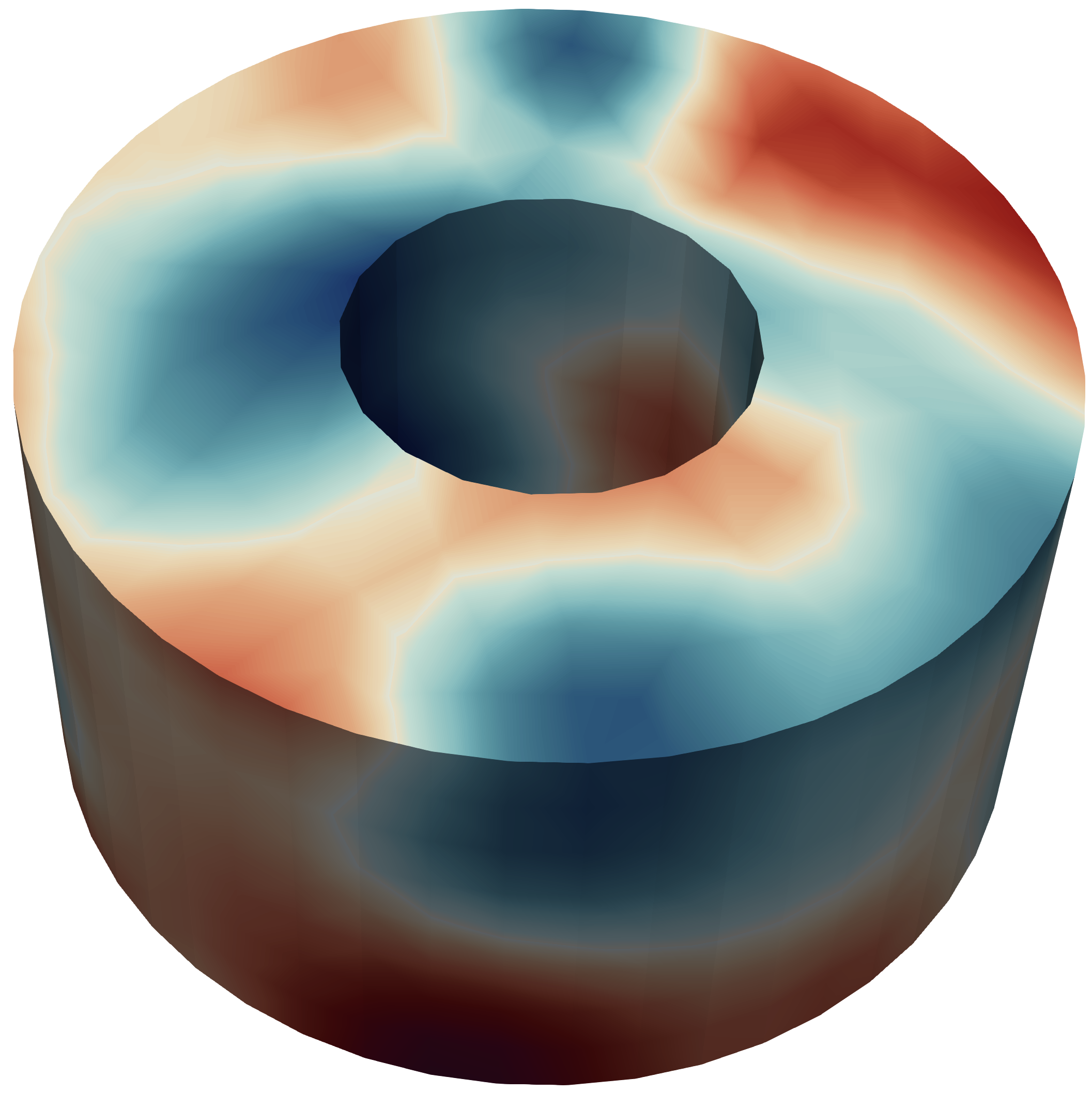}
		      \hspace{0.2cm}
		      \includegraphics[scale=0.0325]{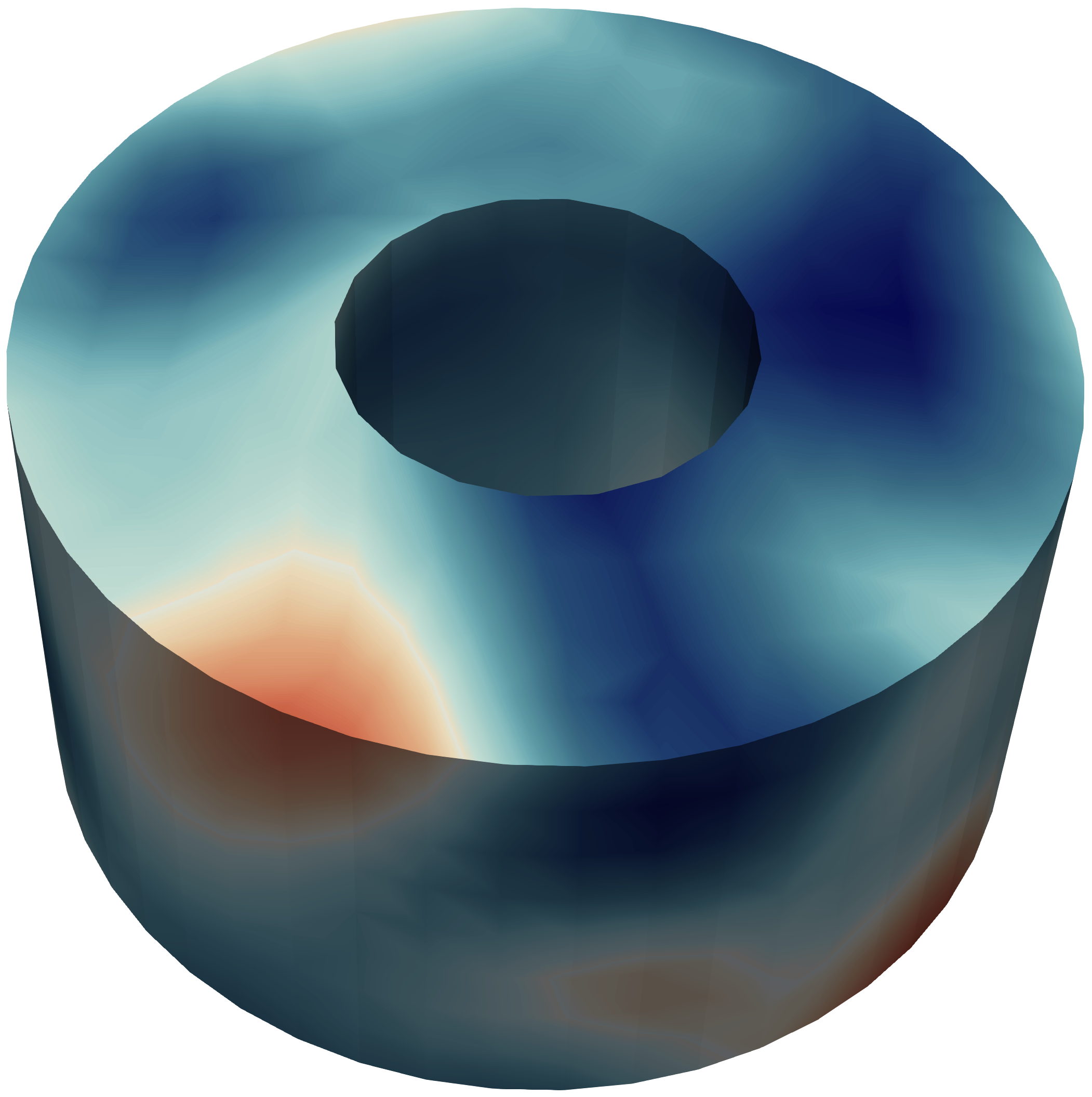}
		      \hspace{0.2cm}
		      \includegraphics[scale=0.0325]{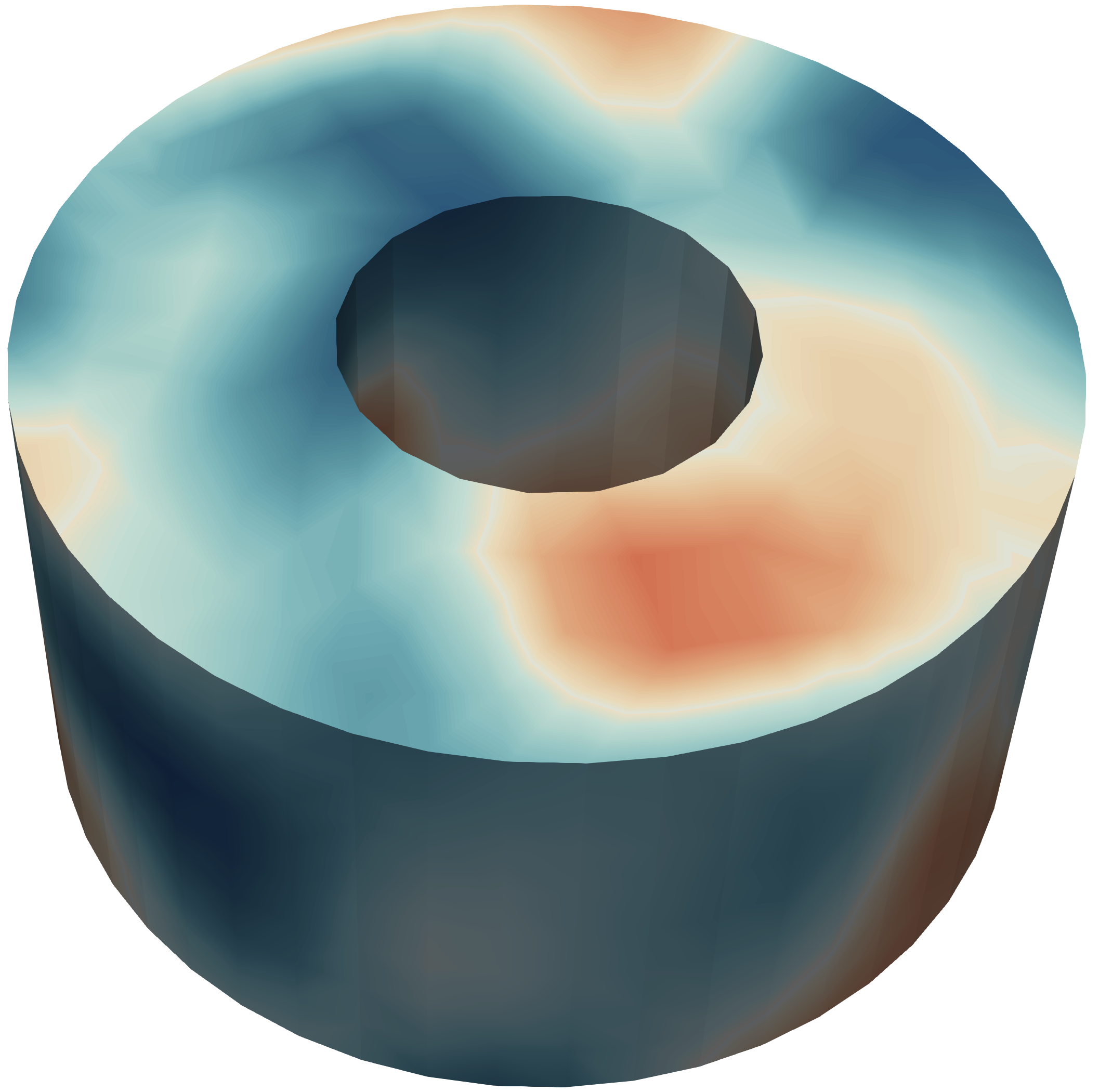}
		      \hspace{0.2cm}
		      \includegraphics[scale=0.0325]{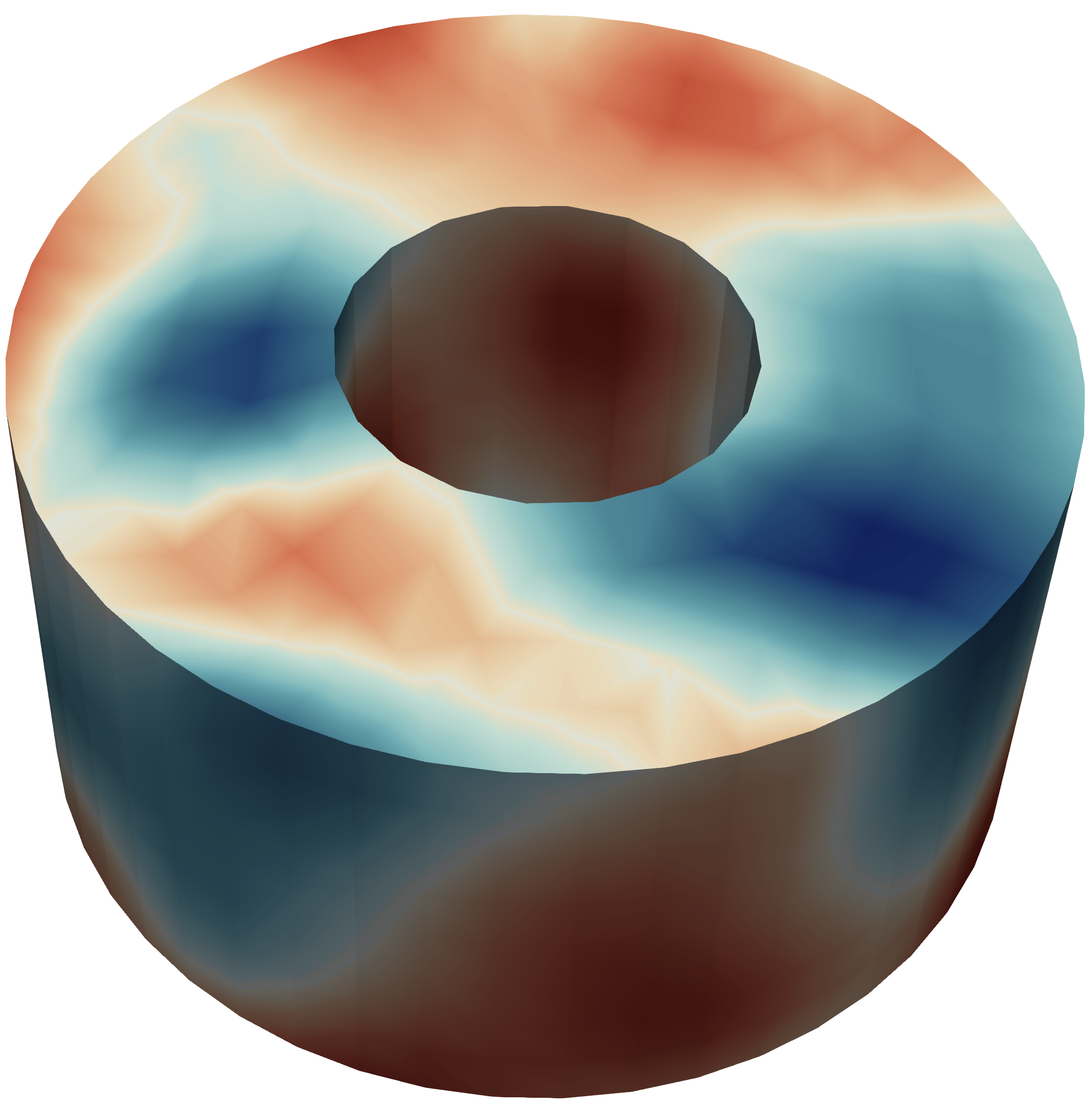}
		      \hspace{0.2cm}
		      \includegraphics[scale=0.065]{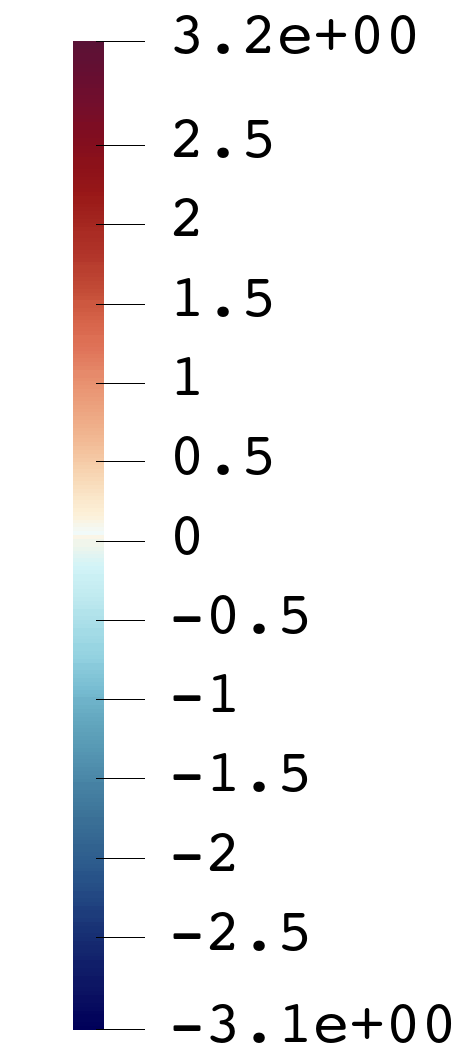}
		      \caption{Example of sampled right-hand sides for creating the Helm dataset.}
		      \label{fig:cylinder}
	      \end{figure}

	\item{\textbf{Helmholtz equation (Helm):}}
	      Let $\Omega$ be the computational domain, defined as an annular cylinder depicted on~\Cref{fig:cylinder}.
	      The parametrized Helmholtz equation is given as
	      \begin{equation}
		      \begin{aligned}
			      - \Delta u(x) - k_{\text{H}}^2 u(x) & = f(x, \boldsymbol{\eta}),  \quad \quad \quad &  & \forall \ x \in \Omega,                      \\
			      u                                   & = 0, \                                        &  & \forall \ x \ \text{on} \ \newtext{\Gamma}\oldtext{\partial\Omega}, \\
                  \newtext{\frac{\partial{u}}{\partial n}}                                   & \newtext{= 0,} \                                        &  & \newtext{\forall \ x \ \text{on} \ \partial \Omega \setminus\Gamma,} \\            
		      \end{aligned}
		      \label{eq:helmholtz}
	      \end{equation}
	      where~$u$ is the solution and \newtext{$\Gamma$}\oldtext{$\partial\Omega$} denotes the top boundary of the cylinder.
	      Moreover, the symbol $k_{\text{H}}$ denotes a constant wave number, which we set to $k_{\text{H}}=3$.
	      We sample the forcing term from Gaussian random fields (GRFs) with mean $\E[f(x, \boldsymbol{\eta})]=0$, and  the covariance is given as
	      \begin{align*}
		      \text{Cov}(f(x,\boldsymbol{\eta}), f(y, \boldsymbol{\eta})) = \sigma^2 \exp \bigg(- \frac{\| x - y \|^2}{2 \ell^2} \bigg),
	      \end{align*} with~$\sigma = 3.0$ and $\ell = 0.1$.
	      To construct the dataset, we discretize~\eqref{eq:helmholtz} using an unstructured mesh with $1,503$ nodes.
	      The sampled right-hand sides are then interpolated onto a $32 \times 32 \times 32$ grid composed of elements within the bounding box $[-1, 1] \times [-1, 1] \times [0, 1]$, generating an input that can be processed by the 3D convolutional layers of the branch network.
          The number of training samples and test samples are $4{,}250$ and $750$, respectively.
	      We have made this new dataset publicly accessible on the ZENODO repository~\cite{zenodo_helm}.

	\item{\textbf{Advection equation (Adv):}}
	      Following~\cite{lu2022comprehensive}, we consider a wave advection equation with periodic boundary conditions, defined at the domain $\Omega := [0,1]$ in space and $[0,1]$ in time.
	      The problem is parameterized with respect to different initial conditions, and it is given as
	      \begin{equation*}
		      \begin{aligned}
			      \frac{\partial u}{\partial t} + \frac{\partial u}{\partial x} & = 0,  \quad \quad \quad         &  & \forall \ (x, t) \in \Omega \times [0, 1], \\
			      u                                                             & = u_0(x, \boldsymbol{\eta}), \  &  & \forall \ (x, t) \in  \Omega \times \{0\}. \\
		      \end{aligned}
	      \end{equation*}
	      The DON is trained to learn a map from the initial condition to the solution of the PDE.
	      To this aim, $u_0$ is parametrized as 
          \begin{equation*}
              u_0(x, \boldsymbol{\eta}) := h_1 1_{\{c_1 - \frac{w}{2}, c_1 + \frac{w}{2}\}} \\ + \sqrt{\max(h_2^2 - a^2(x-c_2)^2, 0)},
          \end{equation*}
          where the parameters $\boldsymbol{\eta} := [a, h_1, h_2, c_1, c_2, w]$ are \newtext{randomly} sampled \newtext{from the uniform distributions $a \sim \pazocal{U}(0.1, 0.2)$, $h_1 \sim \pazocal{U}(0.5, 2.0)$, $h_2 \sim \pazocal{U}(0.5, 2.0)$, $c_1 \sim \pazocal{U}(0.2, 0.3)$, $c_2 \sim \pazocal{U}(0.7, 0.8)$\newtext{, and $w \sim \pazocal{U}(0.1, 0.2)$}, respectively.}\oldtext{randomly from the range $[0.1, 0.2]$, $[0.5, 2.0]$, $[0.5, 2.0]$, $[0.2, 0.3]$, $[0.7, 0.8]$, and $[0.1, 0.2]$, for $a$, $h_1$, $h_2$, $c_1$, $c_2$ and $w$, respectively.}
	      The dataset for training the DON is constructed using a finite difference method with uniform mesh composed of $40 \times 40$ nodes in space-time.
	      The number of training and test samples is $1{,}000$.
	      Note that for this benchmark problem, we utilize the POD-DeepONet architecture~\cite{lu2022comprehensive}, which replaces the trunk network with the POD bases.
	      This allows us to effectivelly demonstrate the capabilities of the proposed layer-wise decomposition for different types of DON architectures.
\end{itemize}

\begin{table}
	\centering
	\caption{The summary of network architectures for DON's benchmark problems\newtext{~(Aniso, Helm and Adv))}.
		The symbol $[ \cdot ]$ indicates the number of channels for \newtext{convolutional neural network~(CNN)} and the width for \newtext{fully-connected neural network~(FNN)}, respectively.
		The symbols {B} and {T} denote the branch network and trunk network, respectively. \newtext{POD-basis is computed by principal component analysis and fixed.}}
	\label{tab:don_architecture}
    \small
	\begin{tabular}{|l||c|r|r|}
		\hline
		{Example}              & \multicolumn{2}{c|}{{Layers}} & {Act.}                                                      \\ \hline \hline
		\multirow{2}{*}{Aniso} & {B}                           & FNN[2, 256, 256, 256, 256, 256, 128]                 & ReLU \\ \cline{2-4}
		                       & {T}                           & FNN[2, 256, 256, 256, 128]                           & Tanh \\ \hline \hline
		\multirow{2}{*}{Helm}  & {B}                           & CNN[1, 40, 60, 100, 180] + FNN[180, 80, 80, 80, 128] & ReLU \\ \cline{2-4}
		                       & {T}                           & FNN[3, 80, 80, 80, 128]                              & Tanh \\ \hline \hline
		\multirow{2}{*}{Adv}   & {B}                           & FNN[40, 64, 64, 64, 64, 64, 64, 64, 32]              & ReLU \\ \cline{2-4}
		                       & {T}                           & POD-basis                                            & -    \\ \hline
	\end{tabular}
\end{table}

\subsection{Implementation details}
\label{sec:impl_details}
All benchmark problems and considered optimizers are implemented on top of the library PyTorch~\cite{paszke2019pytorch}.
The code used for generating the numerical results is part of the open-source library~\texttt{DistTraiNN}~\cite{distrainngit}.
The description of DNN architectures is summarized in~\Cref{tab:pinn_architectures} and~\Cref{tab:don_architecture} for PINNs and DONs, respectively.
All networks are initialized using Xavier initialization strategy~\cite{glorot2010understanding}.
For PINNs, the boundary conditions are imposed by employing the length factor function $\ell$ given as $\ell(x_1, x_2) = x_1(1-x_1)x_2(1-x_2)$ and $\ell(t,x) = t(x+1)(x-1)$  for DA problem and Burg/AC problems, respectively.
Moreover, we employ an adaptive activation function~\cite{jagtap2020adaptive} and a skip-connection~\cite{he2016deep} to prevent gradient vanishing during training and improve the inference performance of the neural network.

As a baseline, we employ the LBFGS optimizer and single-level layer-wise preconditioned LBFGS (SL-LBFGS).
The proposed two-level preconditioner is also tested by preconditioning the LBFGS method, giving rise to the TL-LBFGS optimizer.
All LBFGS variants utilize the limited-memory implementation with the memory size $m=3$.
Moreover, we utilize momentum with parameter $0.9$ and the cubic backtracking line-search with strong Wolfe conditions~\cite[Algorithm A6.3.1, pages 325-327]{dennis1996numerical} to obtain all step-sizes.
For the SL/TL-LBFGS optimizers,  the global and local optimization problems outlined in~\Cref{alg:tlspqn},  are also solved using the LBFGS, which is restarted at every epoch.

Compared to the standard LBFGS algorithm, a single iteration/epoch of the TL-LBFGS method is computationally more expensive.
To analyze the iteration cost,~\Cref{tab:computational_cost} summarizes the gradient evaluations cost $(g_e)$, the update cost (UC) and memory cost (MC)  of the LBFGS and TL-LBFGS optimizers.
As we can see, the LBFGS optimizer requires one gradient evaluation per iteration, i.e.,~$ g_{e} = 1$.
However, UC depends on the size of stored secant pairs ($m$) used for Hessian approximation.
It can be approximately estimated as $4mn$ flops when the compact matrix representation is employed; see~\cite[Sections 3 and 4]{byrd1994representations} for details.
By incorporating the cost associated with evaluating the momentum, the UC equals to~$n+4mn$ flops.
In the perspective of MC, the LBFGS optimizer requires storing two matrices of size $m \times n$ as well as a momentum of length $n$, which takes up memory of size $n+2mn$.

The computational cost of the TL-LBFGS optimizer is associated with global, local and coarse-level training phases.
The global computational cost is exactly same as for the standard LBFGS optimizer.
The gradient evaluation cost for local and coarse-level parts are $k_{s} g_{e_{s}}$ and $k_{0} g_{e_{0}}$, where $(g_{e_{s}})$ and $(g_{e_{0}})$ denote the cost of gradient evaluation of local and coarse-level optimizers, respectively.
Since the number of coarse trainable parameters is the number of subdomain $N_{sd}$, the cost of coarse-level evaluation $g_{e_{0}}$ is generally less than that of local evaluation $g_{e_{s}}$.
In terms of UC, let $\text{UC}_{s}$ and $\text{UC}_{0}$ be the update cost required by the local and coarse optimizers.
Since the UC for $s^{th}$ subnetwork is scaled proportionally to the number of local trainable parameters $n_{s}$ and the number of local iterations $k_{s}$, the local update cost is~$(\sfrac{\tilde{n}_s}{n} )k_s  \text{UC}_{s}$.
Similarly, the coarse-level update cost is~$(\sfrac{n_{0}}{n} ) k_{0}  \text{UC}_{0}$.
Similarly, the coarse-level update cost is $(\sfrac{n_{0}}{n}) k_{0} \,\text{UC}_{0}$.  
Note that the MC of TL-LBFGS can be derived in the same way as the UC, i.e., the local memory cost and the coarse-level memory cost are $(\sfrac{\tilde{n}_s}{n}) \,\text{MC}_{s}$ and $(\sfrac{n_{0}}{n}) \,\text{MC}_{0}$, respectively.  
\newtext{Moreover, note that in the case of the single-level layer-wise preconditioner for L-BFGS, termed SL-LBFGS, the computational cost can be estimated by simply setting $k_{0}=0$.}

\newtext{\begin{remark}
Although the network is decomposed into $N_{sd}$ subnetworks, the application of $\pazocal{M}_{SL/TL}$ requires a full forward/backward pass through the network.
This is reflected in the estimate for the number of gradient evaluations ($g_e$), which is multiplied by the number of local and coarse-level steps, i.e., by~$k_s$ and $k_0$, respectively.
However, the search directions~$\{\dv_{s}\}_{s=1}^{N_{sd}}$ can be computed in parallel, requiring storage only of local gradients and Hessian approximations, as reflected in the estimates for UC and MC.
\end{remark}}


\begin{table}
	\centering
	\caption{Estimates for the number of the gradient evaluations ($ g_{\text{e}}$), update cost (UC), and memory requirements (MC) per iteration for the LBFGS and the TL-LBFGS optimizers.
		The bold symbol denotes the iteration cost associated with incorporating the coarse-level training, i.e., additional cost compared to SL-LBFGS.}
	\label{tab:computational_cost}
	\begin{tabular}{|l|l|l|l|}
		\hline
		Method                    & $ g_{\text{e}}$                   & UC                                                      & MC                                               \\  \hline \hline
		LBFGS                     & 1                                 & $n+ 4mn$                                                & $n+2mn$                                          \\  \hline
		\multirow{3}{*}{TL-LBFGS} & $2 +$                             & $2n+4mn+ $                                              & $2n+2mn+$                                        \\
		                          & $k_{s}( g_{\text{e}_{s}}) $      & $ (\sfrac{\tilde{n}_s}{n} )k_s\text{UC}_{s} +$          & $ (\sfrac{\tilde{n}_s}{n} )\text{{MC}}_{s} + $   \\
		                          & $\bf{ k_{0}(g_{\text{e}_{0}})  }$ & $\bf{(\sfrac{n_{0}}{n} ) k_{0}\text{\textbf{UC}}_{0}} $ & $\bf{(\sfrac{n_{0}}{n} )\text{\textbf{MC}}_{0}}$ \\ \hline
	\end{tabular}
\end{table}

\section{Numerical results}
\label{sec:results}
In this section, we investigate the performance of the proposed layer-wise preconditioner. 
\newtext{Since it has been demonstrated in~\cite{kopanivcakova2024enhancing, kiyani2025optimizer} that, for training SciML models, L-BFGS is clearly superior to Adam and other first-order methods, we focus our numerical results on the performance of the two-level layer-wise preconditioner for L-BFGS, giving rise to TL-LBFGS training algorithm.}
To this aim, we monitor the value of the loss function $\L$ and the relative $L^{2}$ error, given as
\begin{align*}
	\pazocal{E}_{\text{rel}}(u_{\thetav}, u^*) = \frac{ \| u_{\thetav} - u^* \|_{L^2(\Omega)} }{ \| u_{\thetav} \|_{L^2(\Omega)}},
\end{align*}
where $u_{\thetav}$ is the neural network solution and $u^{\ast}$ denotes the exact or the high-fidelity finite element solution.

\subsection{Impact of the overlap on the convergence of SL-LBFGS}
First, we investigate the impact of overlap on the convergence of the SL-LBFGS optimizer.
During this experiment, we fix the number of local iterations at $k_{s}=50$ and consider a varying number of subnetworks $N_{sd} \in \{2, 4, 8\}$ and the size of the overlapping regions $\delta \in \{0, 1, 2\}$.
\Cref{fig:comparison_overlap} illustrates the obtained results for Burg (PINN) and the Adv (DON) equations.
Note, the layer-wise synchronization strategy proposed in~\Cref{sec:recomination_strategy} is not utilized.
Instead, we average the parameters of the layers corresponding to the overlap, an approach traditionally used in numerical analysis.
As we can see from the obtained results, in this particular case, increasing the size of the overlap $\delta$ does not improve convergence, and in some cases, it even causes its deterioration.

\input{chapters/plots/06_perf_wrt_overlap.tex}
\input{chapters/plots/06_perf_wrt_subdomain_sync}

For comparison,~\Cref{fig:comparison_subdomain_sync} demonstrates the obtained results by incorporating the novel layer-wise synchronization strategy.
In contrast to the results reported in~\Cref{fig:comparison_overlap}, we can observe that using the layer-wise synchronization strategy allows for an enhanced convergence as the size of the overlapping region $\delta$ increases.
Here, we, however, point out that as $\delta$ increases, the computational cost of solving the local minimization problems also increases.
Consequently, selecting $\delta=1$ provides the best trade-off between enhanced convergence and the low computational cost.

\subsection{Impact of coarse-level training (SL-LBFGS vs. TL-LBFGS)}
As a next step, we investigate the impact of coarse-level training on the convergence of the proposed preconditioner.
In this experiment, we fix the size of the overlap to $\delta = 1$ and the number of local iterations to $k_{s}=50$.
At the same time, we vary the number of coarse-level iterations $k_{0}$ and the number of subdomains $N_{sd}$, such that $k_{0} \in \{0, 25, 50\}$ and $N_{sd} \in \{2, 4, 8\}$.
\Cref{fig:comparison_coarse} illustrates the obtained results.
As we can see, the loss function decays more rapidly for a sufficient number of coarse iterations $k_{0}$, suggesting that incorporating coarse-level training significantly improves the convergence rate of the proposed TL-LBFGS method.
However, because the coarse-level training process is sequential and cannot be parallelized, an increased number of coarse-level iterations increases the computational cost.
\newtext{Moreover, we would like to remark that the coarse-level model is guaranteed to approximate the fine-level model only in the local neighborhood of the current fine-level iterate and is therefore inherently local.  
Consequently, performing too many $k_{0}$ nonlinear steps on the coarse level may cause the coarse-level model to no longer be consistent with the fine-level model, which can in turn degrade overall convergence. }
To balance computational cost with improvements in convergence speed, we keep $k_{0}$ fairly low for all subsequent numerical experiments, i.e., between $25$ and $50$.

\input{chapters/plots/06_perf_wrt_coarse}

\subsection{Convergence, accuracy and computational cost of LBFGS and SL/TL-LBFGS optimizers}
In this section, we compare the performance of the proposed TL-LBFGS optimizer with the standard LBFGS optimizer.
The comparison is performed with respect to relative error~$\pazocal{E}_{\text{rel}}$, the number of epochs, $g_{e}$ and UC, as discussed in~\Cref{sec:impl_details}.
\newtext{Our first experiment illustrates the impact of preconditioning on the overall performance of training.
We compare the performance of the two-level additive Schwarz method  as a standalone solver (TL-SA) and as a preconditioner for LBFGS (TL-LBFGS). 
The obtained results clearly indicate that the preconditioned variant is significantly more efficient.
These observations are consistent with findings reported in the literature.
For example, V. Dolean et al.~\cite{dolean2016nonlinear} and X. Cai et al.~\cite{cai2002non} demonstrated that the nonlinear two-level additive Schwarz method alone is insufficient to achieve fast convergence.
To overcome this limitation, they used the two-level additive Schwarz to precondition the Newton method.
In view of the superiority of the preconditioning strategy, all remaining experiments employ the proposed two-level additive Schwarz method as a preconditioner for L-BFGS.}

\input{chapters/plots/03_perf_standalone.tex}

\Cref{fig:comparison_pin} and~\Cref{fig:comparison_don} compare the performance of the LBFGS and TL-LBFGS for PINN's and DON's benchmark problems.
As we can see, TL-LBFGS always provides the models with lower~$\pazocal{E}_{\text{rel}}$ than the LBFGS optimizer.
For example, for the DA example with $N_{sd}=5$, the relative error~$\pazocal{E}_{\text{rel}}$ almost immediately decreases to $10^{-2}$, while LBFGS optimizer fails to train the model to that accuracy.
Moreover, we can also observe that for all benchmark problems, the performance of TL-LBFGS improves as the number of subdomains $N_{sd}$ increases.

\input{chapters/plots/06_PIN_BFGS_SPQN.tex}
\input{chapters/plots/06_DON_BFGS_SPQN.tex}

As a next step, we compare the performance of TL-LBFGS with respect to its single-level counterpart, SL-LBFGS.
We configure SL-LBFGS as proposed in~\cite{kopanivcakova2024enhancing}, i.e., without overlap,  without layer-wise synchronization strategy and without the coarse-level training.
As we can see from~\Cref{tab:rel_error}, except in the case of Aniso problem, the TL-LBFGS method consistently outperforms the SL-LBFGS method by providing more accurate DNN models.
While the observed improvements compared to the LBFGS method are around one order of magnitude, they are not as big as those compared to the SL-LBFGS method.
However, in some cases, such as the Burg benchmark problem, where SL-LBFGS did not perform very well~\footnote{
For Burg problem, the SL-LBFGS method starts stagnating after $6{,}000$ epochs at  $\pazocal{E}_{\text{rel}} = 2 \times 10^{-4}$.
Therefore, its overall performance to achieve the same $\pazocal{E}_{\text{rel}}$ as the LBFGS method is poor, despite being very efficient in the initial training phase.}, TL-LBFGS yields about four times more accurate results.


\begin{table}
	\centering
	\caption{Left: The best relative error $\pazocal{E}_{\text{rel}}$ which is achieved by LBFGS, and SL/TL-LBFGS optimizers for PINN's~\newtext{(Allen-Cahn equation~(AC), Diffusion-Advection~(DA), and Burgers'~(Burg))}
	and DON's~\newtext{(Anisotropic Poisson~(Aniso), Advection~(Adv), and Helmholtz~(Helm))} benchmark problem.
	The results are obtained over 5 independent runs.
	Right: An overview of hyperparameters chosen for testing the SL/TL-LBFGS methods. 
	The symbols $N_{sd}$, $k_{s}$, and $k_{0}$ stand for the number of subnetworks,
	the number of iterations for the local problems, and the number of iterations for the coarse problems, respectively.
	\newtext{Note that $N_{sd}$ takes two different values, reflecting the use of two distinct decompositions.  
For example, $2,3$ indicate that $N_{sd}=2$ and $3$ was used.}}

	\label{tab:rel_error}
	\small
	\begin{tabular}{ | c || c | c | c | c | r | r | r |}
		\cline{1-4} \cline{6-8}
		\multicolumn{1}{|c||}{ \multirow{2}{*}{Ex.}} & \multicolumn{3}{c|}{$\pazocal{E}_{\text{rel}}$} &                            & \multicolumn{3}{c|}{Hyperparameters}                                   \\ \cline{2-4} \cline{6-8}
		\multicolumn{1}{|c||}{  }                    & LBFGS                                           & SL-LBFGS                   & TL-LBFGS                             &  & $N_{sd}$ & $k_{s}$ & $k_{0}$ \\ \cline{2-4} \cline{6-8} \cline{1-4} \cline{6-8}
		AC                                           & $3.6 \times 10^{-3}$                            & $2.4 \times 10^{-4}$       & ${\bf 1.6 \times 10^{-4}}$           &  & 2, 6     & 50      & 50      \\ \cline{1-4} \cline{6-8}
		DA                                           & $5.5 \times 10^{-1}$                            & $2.2 \times 10^{-2}$       & ${\bf 1.8 \times 10^{-2}}$           &  & 2, 5     & 50      & 25      \\ \cline{1-4} \cline{6-8}
		Burg                                         & $1.5 \times 10^{-4}$                            & $1.2 \times 10^{-4}$       & ${\bf 2.9 \times 10^{-5}}$           &  & 2, 4     & 50      & 25      \\ \cline{1-4} \cline{6-8}
		Aniso                                        & $1.7 \times 10^{-2}$                            & ${\bf 2.5 \times 10^{-3}}$ & $2.9 \times 10^{-3}$                 &  & 2, 5     & 50      & 50      \\ \cline{1-4} \cline{6-8}
		Adv                                          & $3.0 \times 10^{-2}$                            & $1.9 \times 10^{-2}$       & $\bf{1.1 \times 10^{-2}}$            &  & 2, 4     & 50      & 50      \\ \cline{1-4} \cline{6-8}
		Helm                                         & $7.3 \times 10^{-2}$                            & $2.2 \times 10^{-2}$       & $\bf{1.7 \times 10^{-2}}$            &  & 2, 6     & 25      & 25      \\ \cline{1-4} \cline{6-8}
	\end{tabular}
\end{table}


In the end, we compare $g_{e}$ and UC required to achieve the same accuracy as the LBFGS optimizer.
Tables~\ref{tab:cost_comparison_g} and~\ref{tab:cost_comparison_UC} report the obtained results in terms of speedup factor~$\pazocal{S}$, defined as $\pazocal{S}^A_B = \sfrac{C_A}{C_B}$, where $C_A$ and $C_B$ denote the $g_{e}$ or UC required by the optimizer A and by the optimizer B, respectively.
As we can see, for PINN's problems, which are significantly more ill-conditioned than DON's problems, utilizing the TL-LBFGS optimizer reduces $g_{e}$ and UC by factors of $1.5-1,099$ and $3.5-12,121$ compared to the LBFGS method, respectively.
Moreover, compared to SL-LBFGS, we achieve a speedup of $2.9-7.3$ and $1.6-5.1$ for $g_{e}$ and UC, respectively.
This effectively highlights the benefits of incorporating the proposed coarse-level training step, overlap, and layer-wise synchronization into the training process.

For DON's benchmark problems, we observe lower speedup factors than for the PINNs.
In particular, compared to LBFGS, we obtain speedup factors of $12.8-18.3$ for $g_{e}$ and $33.2-45.5$ for UC, respectively.
Compared to  SL-LBFGS, we notice speedup factors of $0.64-1.7$ for $g_{e}$ and $0.3-1.3$ for UC, respectively.
The most unsatisfactory result is obtained for the Aniso example, where we even observe an increase in the computational cost.
We attribute this behavior to the fact that this problem is the least ill-conditioned of all our benchmark problems, see~\ref{sec:appendix}.
As a consequence, constructing a two-level preconditioner significantly increases the computational cost, while the convergence speed improves only slightly.

\begin{table}
	\centering
	\caption{The number of gradient evaluation~($g_{e}$) required by LBFGS and SL/TL-LBFGS optimizers for achieving the same relative error $\pazocal{E}_{\text{rel}}$ as obtained by the LBFGS optimizer \newtext{for PINN's~(Allen-Cahn~(AC), Diffusion-Advection~(DA), and Burgers'~(Burg))
	and DON's~(Anisotropic Poisson~(Aniso), Advection~(Adv), and Helmholtz~(Helm)) benchmark problem}.
	The results are obtained over 5 independent runs. The speedup factor~$\pazocal{S}$ is defined as $\pazocal{S}^A_B = \sfrac{g_{e}^{A}}{g_{e}^{B}}$, where $g_e^{A}$ and $g_e^{B}$ denote the $g_{e}$ required by the optimizer A and by the optimizer B.}
	\label{tab:cost_comparison_g}
	\small
	\begin{tabular}{| c | c || r || r | r || r | r | r |}
		\hline
		\multicolumn{1}{|c|}{ \multirow{2}{*}{Ex.} } & \multicolumn{1}{c||}{$\pazocal{E}_{\text{rel}}$} & \multicolumn{6}{c|}{$g_{e}$}                                                                                                                                                                                                                                                   \\ \cline{3-8}
		\multicolumn{1}{|c|}{}                       & \multicolumn{1}{c||}{(LBFGS)}                    & \multicolumn{1}{c||}{LBFGS}  & \multicolumn{1}{c|}{SL} & \multicolumn{1}{c||}{$\pazocal{S}_{\text{LBFGS}}^{\text{SL}}$} & \multicolumn{1}{c|}{TL} & \multicolumn{1}{c|}{$\pazocal{S}_{\text{LBFGS}}^{\text{TL}}$} & \multicolumn{1}{c|}{$\pazocal{S}_{\text{SL}}^{\text{TL}}$} \\ \hline \hline
		AC                                           & $3.6 \times 10^{-3}$                             & $1{,}090{,}134$              & $49{,}296$              & $22.11$                                                        & $\boldsymbol{16{,}786}$ & $64.9$                                                        & $2.9$                                                      \\ \hline
		DA                                           & $5.5 \times 10^{-1}$                             & $497{,}000$                  & $1{,}560$               & $318.5$                                                        & $\boldsymbol{452}$      & $1{,}099.5$                                                   & $3.5 $                                                     \\  \hline
		Burg                                         & $1.5 \times 10^{-4}$                             & $99{,}900$                   & $553{,}332$             & 0.2                                                            & $\boldsymbol{75{,}014}$ & 1.5                                                           & 7.3                                                        \\ \hline \hline

		Aniso                                        & $1.7 \times 10^{-2}$                             & $130{,}623$                  & $\boldsymbol{4{,}784}$  & $27.3$                                                         & $7{,}392$               & $17.7$                                                        & $0.64$                                                     \\ \hline
		Adv                                          & $3.0 \times 10^{-2}$                             & $36{,}338$                   & $5{,}044$               & $7.2$                                                          & $\boldsymbol{2{,}838}$  & $12.8$                                                        & $1.7$                                                      \\ \hline
		Helm                                         & $7.3 \times 10^{-2}$                             & $29{,}965$                   & $2{,}673$               & $11.2$                                                         & $\boldsymbol{1{,}638}$  & $18.3$                                                        & $1.6$                                                      \\ \hline
	\end{tabular}
\end{table}

\begin{table}
	\centering
	\caption{The update cost~(UC) required by LBFGS and SL/TL-LBFGS optimizers for achieving the same relative error $\pazocal{E}_{\text{rel}}$ as obtained by the LBFGS optimizer \newtext{for PINN's~(Allen-Cahn~(AC), Diffusion-Advection~(DA), and Burgers'~(Burg))
	and DON's~(Anisotropic Poisson~(Aniso), Advection~(Adv), and Helmholtz~(Helm)) benchmark problem}.
	The results are obtained over 5 independent runs. The speedup factor~$\pazocal{S}$ is defined as $\pazocal{S}^A_B = \sfrac{\text{UC}^{A}}{\text{UC}^{B}}$, where $\text{UC}^{A}$ and $\text{UC}^{B}$ denote the UC required by the optimizer A and by the optimizer B.}
	\label{tab:cost_comparison_UC}
	\footnotesize
	\begin{tabular}{| c | c || r || r | r || r | r | r |}
		\hline
		\multicolumn{1}{|c|}{ \multirow{2}{*}{Ex.} } & \multicolumn{1}{c||}{$\pazocal{E}_{\text{rel}}$} & \multicolumn{6}{c|}{UC}                                                                                                                                                                                                                                                        \\ \cline{3-8}
		\multicolumn{1}{|c|}{}                       & \multicolumn{1}{c||}{(LBFGS)}                    & \multicolumn{1}{c||}{LBFGS} & \multicolumn{1}{c|}{SL} & \multicolumn{1}{c||}{$\pazocal{S}_{\text{LBFGS}}^{\text{SL}}$} & \multicolumn{1}{c|}{TL}  & \multicolumn{1}{c|}{$\pazocal{S}_{\text{LBFGS}}^{\text{TL}}$} & \multicolumn{1}{c|}{$\pazocal{S}_{\text{SL}}^{\text{TL}}$} \\ \hline \hline
		AC                                           & $3.6 \times 10^{-3}$                             & $14{,}171{,}742$            & $167{,}322$             & $84.7$                                                         & $\boldsymbol{106{,}580}$ & $132.9$                                                       & $1.6$                                                      \\ \cline{1-8}
		DA                                           & $5.5 \times 10^{-1}$                             & $6{,}461{,}000$             & $4{,}154$               & $1555.4$                                                       & $\boldsymbol{533}$       & $12,121.9 $                                                   & $2.9$                                                      \\  \cline{1-8}
		Burg                                         & $1.5 \times 10^{-4}$                             & $1{,}298{,}700$             & $1{,}878{,}137$         & $0.69 $                                                        & $\boldsymbol{365{,}496}$ & $3.5 $                                                        & $5.1 $                                                     \\ \hline \hline

		Aniso                                        & $1.7 \times 10^{-2}$                             & $1{,}698{,}099$             & $\boldsymbol{13{,}248}$ & $128.2 $                                                       & $51{,}168$               & $33.2$                                                        & $0.3$                                                      \\ \cline{1-8}
		Adv                                          & $3.0 \times 10^{-2}$                             & $472{,}394$                 & $16{,}072$              & $ 29.4$                                                        & $\boldsymbol{14{,}130}$  & $33.4$                                                        & $1.1$                                                      \\ \cline{1-8}
		Helm                                         & $7.3 \times 10^{-2}$                             & $389{,}545$                 & $11{,}400$              & $34.2$                                                         & $\boldsymbol{8{,}560}$   & $45.5$                                                        & $1.3$                                                      \\ \hline
	\end{tabular}
\end{table}

\section{Discussion}
\newtext{Machine learning (ML)-based solvers for PDEs are rapidly gaining attention; however, their role remains distinct from that of classical numerical methods.
Established approaches such as finite element and finite volume methods, when preconditioned by multigrid or domain decomposition techniques, continue to be the most efficient and reliable tools for solving large-scale, three-dimensional problems on unstructured geometries.
At present, machine learning-based solvers should be regarded primarily as complementary methods, with particular strengths in high-dimensional problems, inverse problems, model discovery, and surrogate modeling for parametric problems.}

\newtext{A major drawback of current machine learning-based approaches is the absence of rigorous error analysis and convergence guarantees, which makes it challenging to assess their reliability compared to classical methods.
Their performance is also highly sensitive to heuristic design choices, including the architecture of the network, the selection of hyperparameters, and the training strategy.
From a computational standpoint, classical methods can often solve benchmark problems in fractions of a second, whereas training machine learning-based solvers typically require significant time and computational resources.
As demonstrated in the manuscript, the development of novel, large-scale, training strategies can significantly reduce the associated computational burden.
Moreover, once trained, neural networks enable extremely fast evaluations across multiple parameter configurations, which is particularly appealing for applications such as shape optimization, design under uncertainty, or real-time simulation.}

\newtext{In summary, although machine learning-based solvers are not yet ready to fully replace state-of-the-art numerical methods, they represent a promising complementary paradigm.  
A key open challenge is the development of novel training strategies that combine improved efficiency with enhanced model accuracy, ultimately enabling machine learning-based solvers to evolve into practical and reliable tools for large-scale scientific computing.}

\section{Conclusion}
\label{sec:conclusion}
In this paper, we proposed a novel two-level overlapping additive Schwarz preconditioner for enhanced training of scientific machine learning applications. 
Motivated by traditional DD methods, the network parameters were decomposed into overlapping groups. 
The parameters of each group were then trained simultaneously, i.e.,~in parallel. 
To achieve fast convergence, we introduced a subdomain-wise synchronization strategy and incorporated a coarse-level training step. 
The performance of the devised two-level preconditioner was investigated using six benchmark problems, encompassing PINNs and DONs. 
The obtained results demonstrated that the proposed preconditioned optimizer yields significantly more accurate models while also reducing the training cost compared to the state-of-the-art LBFGS method and the existent single-level layer-wise preconditioner.

In the future, we plan to combine the proposed two-level layer-wise preconditioning framework with approaches that consider the decomposition of the computational domain/data space, such as those proposed in~\cite{jagtap2021extended,dolean2022finite}. 
Moreover, it would be interesting to extend the proposed training algorithm into stochastic settings, which would enable its applicability beyond scientific machine learning applications. 
Furthermore, developing strategies for enhanced parallel implementation, such as load balancing and optimization of the subnetwork's loss/gradient evaluations, could further improve the performance.

\section*{Acknowledgements}

Y.L. was supported in part by Basic Science Research Program through NRF funded by the Ministry of Education~(No. RS2023-00247199).  
A.K. was supported by the Swiss National Science Foundation~(SNSF) under the project ``\textit{Multilevel training of DeepONets – multiscale and multiphysics applications}'', as well as by the Platform for Advanced Scientific Computing~(PASC) under the project \textit{ExaTrain}. 
Moreover, work of A.K. benefited from the AI Interdisciplinary Institute ANITI, funded by the  France 2030 program under Grant Agreement No. ANR-23-IACL-0002.
G.E.K. is supported by the ONR Vannevar Bush Faculty Fellowship.
We also acknowledge support from the DOE-MMICS SEA-CROGS DE-SC0023191 award and Ansys Inc.

\bibliographystyle{elsarticle-num}
\bibliography{biblio_cleaned}

\begin{thebibliography}{10}
\expandafter\ifx\csname url\endcsname\relax
  \def\url#1{\texttt{#1}}\fi
\expandafter\ifx\csname urlprefix\endcsname\relax\def\urlprefix{URL }\fi
\expandafter\ifx\csname href\endcsname\relax
  \def\href#1#2{#2} \def\path#1{#1}\fi

\bibitem{barron1993universal}
A.~R. Barron, Universal approximation bounds for superpositions of a sigmoidal function, IEEE Transactions on Information Theory 39 (1993) 930--945.

\bibitem{chen1990constructive}
T.~Chen, H.~Chen, R.-w. Liu, A constructive proof of {C}ybenko's approximation theorem, in: Computing Science and Statistics, Springer New York, NY, 1992, pp. 163--168.

\bibitem{siegel2020approximation}
J.~W. Siegel, J.~Xu, Approximation rates for neural networks with general activation functions, Neural Networks 128 (2020) 313--321.

\bibitem{karniadakis2021physics}
G.~E. Karniadakis, I.~G. Kevrekidis, L.~Lu, P.~Perdikaris, S.~Wang, L.~Yang, Physics-informed machine learning, Nature Reviews Physics 3~(6) (2021) 422--440.

\bibitem{raissi2019physics}
M.~Raissi, P.~Perdikaris, G.~E. Karniadakis, Physics-informed neural networks: A deep learning framework for solving forward and inverse problems involving nonlinear partial differential equations, Journal of Computational Physics 378 (2019) 686--707.

\bibitem{sirignano2018dgm}
J.~Sirignano, K.~Spiliopoulos, {DGM}: A deep learning algorithm for solving partial differential equations, Journal of Computational Physics 375 (2018) 1339--1364.

\bibitem{mishra2022estimates}
S.~Mishra, R.~Molinaro, Estimates on the generalization error of physics-informed neural networks for approximating a class of inverse problems for {PDEs}, IMA Journal of Numerical Analysis 42~(2) (2022) 981--1022.

\bibitem{lu2021learning}
L.~Lu, P.~Jin, G.~Pang, Z.~Zhang, G.~E. Karniadakis, Learning nonlinear operators via {D}eep{O}{N}et based on the universal approximation theorem of operators, Nature Machine Intelligence 3~(3) (2021) 218--229.

\bibitem{wang2021learning}
S.~Wang, H.~Wang, P.~Perdikaris, Learning the solution operator of parametric partial differential equations with physics-informed {D}eep{ON}ets, Science Advances 7~(40) (2021) eabi8605.

\bibitem{li2021fourier}
Z.~Li, N.~B. Kovachki, K.~Azizzadenesheli, B.~liu, K.~Bhattacharya, A.~Stuart, A.~Anandkumar, Fourier neural operator for parametric partial differential equations, in: International Conference on Learning Representations, 2021.

\bibitem{kingma2014adam}
D.~P. Kingma, J.~Ba, Adam: A method for stochastic optimization, arXiv preprint arXiv:1412.6980 (2014).

\bibitem{liu1989limited}
D.~C. Liu, J.~Nocedal, On the limited memory {BFGS} method for large scale optimization, Mathematical Programming 45 (1989) 503--528.

\bibitem{Liu2020On}
L.~Liu, H.~Jiang, P.~He, W.~Chen, X.~Liu, J.~Gao, J.~Han, On the variance of the adaptive learning rate and beyond, in: International Conference on Learning Representations, 2020.

\bibitem{loshchilov2018decoupled}
I.~Loshchilov, F.~Hutter, Decoupled weight decay regularization, in: International Conference on Learning Representations, 2019.

\bibitem{Al-Baali03092014}
M.~Al-Baali, E.~Spedicato, F.~Maggioni, Broyden's quasi-newton methods for a nonlinear system of equations and unconstrained optimization: a review and open problems, Optimization Methods and Software 29~(5) (2014) 937--954.

\bibitem{ben2019demystifying}
T.~Ben-Nun, T.~Hoefler, Demystifying parallel and distributed deep learning: An in-depth concurrency analysis, ACM Computing Surveys (CSUR) 52~(4) (2019) 1--43.

\bibitem{dean2012large}
J.~Dean, G.~Corrado, R.~Monga, K.~Chen, M.~Devin, M.~Mao, M.~a. Ranzato, A.~Senior, P.~Tucker, K.~Yang, Q.~Le, A.~Ng, Large scale distributed deep networks, in: Advances in Neural Information Processing Systems, Vol.~25, 2012, pp. 1223--1231.

\bibitem{krizhevsky2014one}
A.~Krizhevsky, One weird trick for parallelizing convolutional neural networks, arXiv preprint arXiv:1404.5997 (2014).

\bibitem{huang2019gpipe}
Y.~Huang, Y.~Cheng, A.~Bapna, O.~Firat, D.~Chen, M.~Chen, H.~Lee, J.~Ngiam, Q.~V. Le, Y.~Wu, et~al., G{P}ipe: Efficient training of giant neural networks using pipeline parallelism, in: Advances in Neural Information Processing Systems, Vol.~32, 2019, pp. 103--112.

\bibitem{narayanan2019pipedream}
D.~Narayanan, A.~Harlap, A.~Phanishayee, V.~Seshadri, N.~R. Devanur, G.~R. Ganger, P.~B. Gibbons, M.~Zaharia, Pipe{D}ream: {G}eneralized pipeline parallelism for dnn training, in: ACM Symposium on Operating Systems Principles, 2019, pp. 1--15.

\bibitem{rudakov2024activations}
M.~Rudakov, A.~Beznosikov, Y.~Kholodov, A.~Gasnikov, Activations and gradients compression for model-parallel training, Doklady Mathematics 108 (2023) S272--S281.

\bibitem{dolean2015introduction}
V.~Dolean, P.~Jolivet, F.~Nataf, An introduction to domain decomposition methods: algorithms, theory, and parallel implementation, SIAM Philadelphia, 2015.

\bibitem{toselli2004domain}
A.~Toselli, O.~Widlund, Domain decomposition methods-algorithms and theory, Springer Berlin, 2004.

\bibitem{lions1988schwarz}
P.-L. Lions, et~al., On the {S}chwarz alternating method. {I}, in: First International Symposium on Domain Decomposition Methods for Partial Differential Equations, Vol.~1, Paris, France, 1988, p.~42.

\bibitem{dryja1989additive}
M.~Dryja, An additive {S}chwarz algorithm for two-and three-dimensional finite element elliptic problems, in: Domain Decomposition Methods, SIAM Philadelphia, 1989, pp. 168--172.

\bibitem{dryja1989towards}
M.~Dryja, O.~B. Widlund, Towards a unified theory of domain decomposition algorithms for elliptic problems, in: Third International Symposium on Domain Decomposition Methods for Partial Differential Equations, SIAM Philadelphia, 1990, pp. 3--21.

\bibitem{gu2022decomposition}
L.~Gu, W.~Zhang, J.~Liu, X.-C. Cai, Decomposition and composition of deep convolutional neural networks and training acceleration via sub-network transfer learning, Electronic Transactions on Numerical Analysis (2022) 157--186.

\bibitem{klawonn2023domain}
A.~Klawonn, M.~Lanser, J.~Weber, A domain decomposition–based {CNN-DNN} architecture for model parallel training applied to image recognition problems, SIAM Journal on Scientific Computing 46~(5) (2024) C557--C582.

\bibitem{gunther2020layer}
S.~G\"{u}nther, L.~Ruthotto, J.~B. Schroder, E.~C. Cyr, N.~R. Gauger, Layer-parallel training of deep residual neural networks, SIAM Journal on Mathematics of Data Science 2~(1) (2020) 1--23.

\bibitem{kopanivcakova2024enhancing}
A.~Kopani{\v{c}}{\'a}kov{\'a}, H.~Kothari, G.~E. Karniadakis, R.~Krause, Enhancing training of physics-informed neural networks using domain decomposition--based preconditioning strategies, SIAM Journal on Scientific Computing (2024) S46--S67.

\bibitem{lee2022parareal}
Y.~Lee, J.~Park, C.-O. Lee, Parareal neural networks emulating a parallel-in-time algorithm, IEEE Transactions on Neural Networks and Learning Systems (2022).

\bibitem{jagtap2021extended}
A.~D. Jagtap, G.~E. Karniadakis, Extended physics-informed neural networks ({XPINNs}): A generalized space-time domain decomposition based deep learning framework for nonlinear partial differential equations, Communications in Computational Physics 28~(5) (2020).

\bibitem{jagtap2020conservative}
A.~D. Jagtap, E.~Kharazmi, G.~E. Karniadakis, Conservative physics-informed neural networks on discrete domains for conservation laws: Applications to forward and inverse problems, Computer Methods in Applied Mechanics and Engineering 365 (2020) 113028.

\bibitem{hu2022augmented}
Z.~Hu, A.~D. Jagtap, G.~E. Karniadakis, K.~Kawaguchi, Augmented physics-informed neural networks ({APINNs}): A gating network-based soft domain decomposition methodology, Engineering Applications of Artificial Intelligence 126 (2023) 107183.

\bibitem{dolean2023multilevel}
V.~Dolean, A.~Heinlein, S.~Mishra, B.~Moseley, Multilevel domain decomposition-based architectures for physics-informed neural networks, SAM Research Report (2023).

\bibitem{LEE2022323}
Y.~Lee, J.~Park, C.-O. Lee, Two-level group convolution, Neural Networks 154 (2022) 323--332.

\bibitem{cai2002non}
X.-C. Cai, D.~E. Keyes, L.~Marcinkowski, Non-linear additive schwarz preconditioners and application in computational fluid dynamics, International Journal for Numerical Methods in Fluids 40~(12) (2002) 1463--1470.

\bibitem{dolean2016nonlinear}
V.~Dolean, M.~J. Gander, W.~Kheriji, F.~Kwok, R.~Masson, Nonlinear preconditioning: How to use a nonlinear {Schwarz} method to precondition {Newton}'s method, SIAM Journal on Scientific Computing 38~(6) (2016) A3357--A3380.

\bibitem{kopanivcakova2023nonlinear}
A.~Kopani{\v{c}}{\'a}kov{\'a}, H.~Kothari, R.~Krause, Nonlinear field-split preconditioners for solving monolithic phase-field models of brittle fracture, Computer Methods in Applied Mechanics and Engineering 403 (2023) 115733.

\bibitem{kothari2022nonlinear_bounds}
H.~Kothari, A.~Kopani\v{c}{\'a}kov{\'a}, R.~Krause, Nonlinear {S}chwarz preconditioning for nonlinear optimization problems with bound constraints, in: Domain Decomposition Methods in Science and Engineering XXVII, Springer Nature Switzerland, 2024, pp. 319--326.

\bibitem{luo2023preconditioned}
L.~Luo, X.-C. Cai, Preconditioned inexact {Newton} with learning capability for nonlinear system of equations, SIAM Journal on Scientific Computing 45~(2) (2023) A849--A871.

\bibitem{chaouqui2022linear}
F.~Chaouqui, M.~J. Gander, P.~M. Kumbhar, T.~Vanzan, Linear and nonlinear substructured restricted additive {S}chwarz iterations and preconditioning, Numerical Algorithms 91~(1) (2022) 81--107.

\bibitem{he2016deep}
K.~He, X.~Zhang, S.~Ren, J.~Sun, Deep residual learning for image recognition, in: Proceedings of the IEEE Conference on Computer Vision and Pattern Recognition, 2016, pp. 770--778.

\bibitem{sukumar2022exact}
N.~Sukumar, A.~Srivastava, Exact imposition of boundary conditions with distance functions in physics-informed deep neural networks, Computer Methods in Applied Mechanics and Engineering 389 (2022) 114333.

\bibitem{wang2022and}
S.~Wang, X.~Yu, P.~Perdikaris, When and why {PINNs} fail to train: A neural tangent kernel perspective, Journal of Computational Physics 449 (2022) 110768.

\bibitem{anagnostopoulos2024residual}
S.~J. Anagnostopoulos, J.~D. Toscano, N.~Stergiopulos, G.~E. Karniadakis, Residual-based attention in physics-informed neural networks, Computer Methods in Applied Mechanics and Engineering 421 (2024) 116805.

\bibitem{lu2021physics}
L.~Lu, R.~Pestourie, W.~Yao, Z.~Wang, F.~Verdugo, S.~G. Johnson, Physics-informed neural networks with hard constraints for inverse design, SIAM Journal on Scientific Computing 43~(6) (2021) B1105--B1132.

\bibitem{lu2022comprehensive}
L.~Lu, X.~Meng, S.~Cai, Z.~Mao, S.~Goswami, Z.~Zhang, G.~E. Karniadakis, A comprehensive and fair comparison of two neural operators (with practical extensions) based on fair data, Computer Methods in Applied Mechanics and Engineering 393 (2022) 114778.

\bibitem{anandkumar2019neural}
A.~Anandkumar, K.~Azizzadenesheli, K.~Bhattacharya, N.~Kovachki, Z.~Li, B.~Liu, A.~Stuart, Neural operator: Graph kernel network for partial differential equations, in: ICLR 2020 Workshop on Integration of Deep Neural Models and Differential Equations, 2020.

\bibitem{wang2021understanding}
S.~Wang, Y.~Teng, P.~Perdikaris, Understanding and mitigating gradient flow pathologies in physics-informed neural networks, SIAM Journal on Scientific Computing 43~(5) (2021) A3055--A3081.

\bibitem{matsokin1985schwarz}
A.~M. Matsokin, S.~V. Nepomnyaschikh, The {S}chwarz alternation method in a subspace, Izvestiya Vysshikh Uchebnykh Zavedenii. Matematika~(10) (1985) 61--66.

\bibitem{bridson2006multipreconditioned}
R.~Bridson, C.~Greif, A multipreconditioned conjugate gradient algorithm, SIAM Journal on Matrix Analysis and Applications 27~(4) (2006) 1056--1068.

\bibitem{greif2017gmres}
C.~Greif, T.~Rees, D.~B. Szyld, {GMRES} with multiple preconditioners, SeMA Journal 74 (2017) 213--231.

\bibitem{brezinski1999multiparameter}
C.~Brezinski, Multiparameter descent methods, Linear Algebra and its Applications 296~(1-3) (1999) 113--141.

\bibitem{kunstner2024searching}
F.~Kunstner, V.~Sanches~Portella, M.~Schmidt, N.~Harvey, Searching for optimal per-coordinate step-sizes with multidimensional backtracking, in: Advances in Neural Information Processing Systems, Vol.~36, 2024, pp. 2725--2767.

\bibitem{gaedke2021multilevel}
L.~Gaedke-Merzh{\"a}user, A.~Kopani{\v{c}}{\'a}kov{\'a}, R.~Krause, Multilevel minimization for deep residual networks, ESAIM: Proceedings and Surveys 71 (2021) 131--144.

\bibitem{gratton2023multilevel}
S.~Gratton, A.~Kopani\v{c}{\'a}kov{\'a}, P.~L. Toint, Multilevel objective-function-free optimization with an application to neural networks training, SIAM Journal on Optimization 33~(4) (2023) 2772--2800.

\bibitem{kopanicakova2022globally}
A.~{Kopani{\v{c}}{\'a}kov{\'a}}, R.~Krause, Globally convergent multilevel training of deep residual networks, SIAM Journal on Scientific Computing (2022) S254--S280.

\bibitem{nocedal2006numerical}
J.~Nocedal, S.~J. Wright, Numerical Optimization, 2nd Edition, Springer, 2006.
\newblock \href {https://doi.org/10.1007/978-0-387-40065-5} {\path{doi:10.1007/978-0-387-40065-5}}.

\bibitem{morokoff1995quasi}
W.~J. Morokoff, R.~E. Caflisch, Quasi-{M}onte {C}arlo integration, Journal of Computational Physics 122 (1995) 218--230.

\bibitem{asmussen2007stochastic}
S.~Asmussen, P.~W. Glynn, Stochastic simulation: algorithms and analysis, Vol.~57, Springer New York, NY, 2007.

\bibitem{FiredrakeUserManual}
D.~A. Ham, P.~H.~J. Kelly, L.~Mitchell, C.~J. Cotter, R.~C. Kirby, K.~Sagiyama, N.~Bouziani, S.~Vorderwuelbecke, T.~J. Gregory, J.~Betteridge, D.~R. Shapero, R.~W. Nixon-Hill, C.~J. Ward, P.~E. Farrell, P.~D. Brubeck, I.~Marsden, T.~H. Gibson, M.~Homolya, T.~Sun, A.~T.~T. McRae, F.~Luporini, A.~Gregory, M.~Lange, S.~W. Funke, F.~Rathgeber, G.-T. Bercea, G.~R. Markall, Firedrake User Manual, Imperial College London and University of Oxford and Baylor University and University of Washington, first edition Edition (5 2023).
\newblock \href {https://doi.org/10.25561/104839} {\path{doi:10.25561/104839}}.

\bibitem{zenodo_aniso}
A.~Kopani{\v c}{\'a}kov{\'a}, \href{https://doi.org/10.5281/zenodo.10909052}{{DON dataset: AnisotropicDiffusion2D}}, Zenodo (2024).
\newblock \href {https://doi.org/https://doi.org/10.5281/zenodo.10909052} {\path{doi:https://doi.org/10.5281/zenodo.10909052}}.
\newline\urlprefix\url{https://doi.org/10.5281/zenodo.10909052}

\bibitem{zenodo_helm}
A.~Kopani{\v c}{\'a}kov{\'a}, \href{https://doi.org/10.5281/zenodo.10904349}{{DON dataset: NonNestedHelmholtz3D}}, Zenodo (2024).
\newblock \href {https://doi.org/https://doi.org/10.5281/zenodo.10904349} {\path{doi:https://doi.org/10.5281/zenodo.10904349}}.
\newline\urlprefix\url{https://doi.org/10.5281/zenodo.10904349}

\bibitem{paszke2019pytorch}
A.~Paszke, S.~Gross, F.~Massa, A.~Lerer, J.~Bradbury, G.~Chanan, T.~Killeen, Z.~Lin, N.~Gimelshein, L.~Antiga, A.~Desmaison, A.~K\"{o}pf, E.~Yang, Z.~DeVito, M.~Raison, A.~Tejani, S.~Chilamkurthy, B.~Steiner, L.~Fang, J.~Bai, S.~Chintala, {PyTorch}: An imperative style, high-performance deep learning library, in: Advances in Neural Information Processing Systems, Vol.~32, 2019, pp. 8026--8037.

\bibitem{distrainngit}
A.~Kopani{\v c}{\'a}kov{\'a}, Y.~Lee, H.~Kothari, \href{https://bitbucket.org/alena\_kopanicakova/disttrainn}{{DistTraiNN: Model parallel framework for distributed training of scientific machine-learning applications. {G}it repository}}, https://bitbucket.org/alena\_kopanicakova/disttrainn (2023).
\newline\urlprefix\url{https://bitbucket.org/alena\_kopanicakova/disttrainn}

\bibitem{glorot2010understanding}
X.~Glorot, Y.~Bengio, Understanding the difficulty of training deep feedforward neural networks, in: Proceedings of the Thirteenth International Conference on Artificial Intelligence and Statistics, Vol.~9, PMLR, 2010, pp. 249--256.

\bibitem{jagtap2020adaptive}
A.~D. Jagtap, K.~Kawaguchi, G.~E. Karniadakis, Adaptive activation functions accelerate convergence in deep and physics-informed neural networks, Journal of Computational Physics 404 (2020) 109136.

\bibitem{dennis1996numerical}
J.~E. Dennis, R.~B. Schnabel, Numerical Methods for Unconstrained Optimization and Nonlinear Equations, SIAM Philadelphia, 1996.

\bibitem{byrd1994representations}
R.~H. Byrd, J.~Nocedal, R.~B. Schnabel, Representations of quasi-{N}ewton matrices and their use in limited memory methods, Mathematical Programming 63~(1-3) (1994) 129--156.

\bibitem{kiyani2025optimizer}
E.~Kiyani, K.~Shukla, J.~F. Urb{\'a}n, J.~Darbon, G.~E. Karniadakis, Optimizing the optimizer for physics-informed neural networks and kolmogorov-arnold networks, Computer Methods in Applied Mechanics and Engineering 446 (2025) 118308.

\bibitem{dolean2022finite}
V.~Dolean, A.~Heinlein, S.~Mishra, B.~Moseley, Finite basis physics-informed neural networks as a {S}chwarz domain decomposition method, in: Domain Decomposition Methods in Science and Engineering XXVII, Springer Nature Switzerland, 2024, pp. 165--172.

\bibitem{de2023operator}
T.~De~Ryck, F.~Bonnet, S.~Mishra, E.~de~B{\'e }zenac, An operator preconditioning perspective on training in physics-informed machine learning, in: The Twelfth International Conference on Learning Representations, 2024.

\bibitem{jacot2018neural}
A.~Jacot, F.~Gabriel, C.~Hongler, {Neural tangent kernel: Convergence and generalization in neural networks}, in: Advances in Neural Information Processing Systems, Vol.~31, 2018, pp. 8580--8589.

\end{thebibliography}

\appendix
\section{Condition number of (layer-wise) neural tangent kernel}
\label{sec:appendix}
\Cref{tab:cond_number} reports the logarithmic value of condition number of the neural tangent kernel~(NTK)~\cite{de2023operator, jacot2018neural}, denoted by log(cond$(\mathbb{A})$) for all benchmark problems considered in~\Cref{sec:examples}.
As we can see, the problems associated with training PINNs are more ill-conditioned than those associated with training DONs.
Moreover, we also see that for PINN problems, the layer-wise condition number, denoted by \(\text{cond}(\mathbb{A}_l)\), varies over a wider range.
This suggests that the local training associated with solving~\eqref{eq:local} requires different adjustments of step-size. Thus, by decoupling the layers, we are able to use larger step-sizes for layers that are better conditioned, which in turn enhances the overall convergence.

\begin{table}[t]
	\centering
	\caption{Average condition number of the NTK for all benchmark problems.}
	\label{tab:cond_number}
	\small
	\begin{tabular}{| r | c | c | c |}
		\hline
		\multirow{2}{*}{{Ex.} } & \multirow{2}{*}{{$\approx$ log(cond$(\mathbb{A})$)}} & \multicolumn{2}{c|}{{Layer-wise $\approx$ log(cond$(\mathbb{A}_{l}))$}}                            \\  \cline{3-4}
		                        &                                                      & \ \ \ \ \ min \ \ \ \ \                                                 & \ \ \ \ \ max  \ \ \ \ \ \\ \hline \hline
		AC                      & $26$                                                 & $21$                                                                    & $27$                     \\    \hline
		DA                      & $23$                                                 & $20$                                                                    & $26$                     \\  \hline
		Burg                    & $33$                                                 & $29$                                                                    & $36$                     \\  \hline  \hline

		Aniso                   & $23$                                                 & $21$                                                                    & $23$                     \\  \hline
		Adv                     & $22$                                                 & $15$                                                                    & $22$                     \\  \hline
		Helm                    & $21$                                                 & $18$                                                                    & $21$                     \\  \hline
	\end{tabular}
\end{table}

\section{Getting started with DistTrainNN}
\newtext{In this section, we briefly describe how to run an TL-LBFGS training algorithms within our open-source ~\texttt{DistTraiNN}~\cite{distrainngit}.
Firstly, the user is required to setup a Python environment with the following depedencies: PyTorch, Pandas, pyDOE, scikit-optimize, scikit-learn.}


\newtext{The scripts used to obtain the numerical results reported in Section~\ref{sec:results} are available in the folder ``tests/SPQN\_tests''.  
An example code for training a PINN to approximate a solution for the Allen--Cahn equation is provided below.
Please, note that the optimizer is fully configurable via command line arguments. }

\begin{lstlisting}[language=bash]

git clone git@bitbucket.org:alena_kopanicakova/disttrainn.git
cd disttrainn

export PYTHONPATH="${PYTHONPATH}:/path/to/disttrainn"

cd tests/TL_SPQN_tests/AllenCahn
python3 test_SPQN_twolevel_additive_serial.py --seed 1 --num_points_x 25 --num_points_y 50 --num_ref 4 --max_epochs 200000 --lr_global 1.0 --lr_local 1.0 --history 3 --width 64 --hiddenlayers_coarse 4 --num_subdomains 3 --num_local_steps 50 --num_coarse_steps 50 --use_layerwise_lr True --overlap_width 1
\end{lstlisting}

\end{document}